\documentclass[leqno]{amsart}
\usepackage{amsfonts,amssymb,amsmath,amsgen,amsthm}
\usepackage{hyperref}
\usepackage{graphicx,subfigure}

\theoremstyle{plain}
\newtheorem{theorem}{Theorem}[section]

\newtheorem{lemme}[theorem]{Lemma}

\newtheorem{proposition}[theorem]{Proposition}

\theoremstyle{remark}
\newtheorem{remark}[theorem]{Remark}

\def\R{{\mathbf R}}
\def\N{{\mathbf N}}

\def\({\left(}
\def\){\right)}
\def\<{\left\langle}
\def\>{\right\rangle}

\def\ge{\geqslant}

\def\Tend#1#2{\mathop{\longrightarrow}\limits_{#1\rightarrow#2}}

\def\d{{\partial}}
\def\eps{\varepsilon}

\def\R{{\mathbb R}}
\def\N{{\mathbb N}}

\def\d{{\partial}}

\def\eps{\varepsilon}
\def\op_#1{\mathrel{\mathop{{\rm op}_{#1}}}}

\def\build#1_#2^#3{\mathrel{
\mathop{\kern 0pt#1}\limits_{#2}^{#3}}}

\def\td_#1,#2{\mathrel{
\mathop{\build\longrightarrow_{#1\rightarrow #2}^{}}}}

\def\lim_#1,#2{\mathrel{
\mathop{\build{\rm lim}_{#1\rightarrow#2}^{}}}}

\def\limsup_#1,#2{\mathrel{
\mathop{\build{\rm limsup}_{#1\rightarrow#2}^{}}}}

\def\liminf_#1,#2{\mathrel{
\mathop{\build{\rm liminf}_{#1\rightarrow#2}^{}}}}

\def\aref#1{(\ref{#1})}
\def\eps{\varepsilon}

\def\op{{\rm op}}

\def\d{{\rm d}}

\numberwithin{equation}{section}

\begin{document}

\title[Wigner measure propagation and Lipschitz singularity]{Wigner measure propagation and conical singularity for general initial data}
\author[C. Fermanian]{Clotilde~Fermanian-Kammerer}
\address[C. Fermanian]{LAMA, UMR CNRS 8050,
Universit\'e Paris EST\\
61, avenue du G\'en\'eral de Gaulle\\
94010 Cr\'eteil Cedex\\ France}
\email{Clotilde.Fermanian@u-pec.fr}
\author[P. G\'erard]{Patrick G\'erard}
\address[P. G\'erard]{Universit\'e  Paris-Sud, Math\'ematiques, Bat. 425, 91405 ORSAY,
FRANCE}
\email{Patrick.Gerard@math.u-psud.fr}
\author[C. Lasser]{Caroline Lasser}
\address[C. Lasser]{Zentrum Mathematik - M3,
Wissenschaftliches Rechnen,
Technische Universit\"at M\"unchen,
85747 Garching bei M\"unchen, GERMANY}
\email{classer@ma.tum.de}

\begin{abstract}
  We study the evolution of Wigner measures of a family of solutions of a Schr\"odinger  equation with a scalar potential displaying a conical singularity. Under a genericity assumption, classical trajectories  exist and are unique, thus the question of the propagation of Wigner measures along these trajectories becomes relevant. We  prove the propagation for general initial data. \end{abstract}
\thanks{}  
\maketitle


\section{Introduction}

We consider  the Schr\"odinger equation 
\begin{equation}\label{eq:1}
\left\{\begin{array}{l}
i\eps\partial_t  \psi^\eps=-\frac{\eps^2}{2}\Delta\psi^\eps +V(x)\psi^\eps, \;\;(t,x)\in \R\times \R^d,\\
\psi^\eps_{|t=0}=\psi^\eps_0.\end{array}\right. 
\end{equation}
where the potential $V(x)$ displays a conical singularity: there exist two scalar-valued functions $w,V_0\in{\mathcal C}^\infty(\R^d,\R)$ and a vector-valued function  $g\in{\mathcal C}^\infty(\R^d,\R^p)$, $0<p\leq d$ such that $g=0$ is a system of equations of  a codimension~$p$ submanifold of $\R^d$ and 
\begin{equation}\label{eq:sing}
\forall x\in\R^d,\;\;V(x)=w(x)|g(x)|+V_0(x).
\end{equation}
We suppose that $V$ satisfies Kato conditions (see~\cite{Kato72}) so that the Schr\"odinger operator $-\frac{\eps^2}{2}\Delta +V(x)$ is essentially  self-adjoint.
Moreover, we are concerned with the effects of conical singularities in the potential; therefore, we assume that $V_0$ and $w$ are smooth. This smoothness assumption can be slightly relaxed as discussed in Remarks~\ref{rem:smooth1} and~\ref{rem:smooth2} below. \\

\noindent Assuming that $(\psi^\eps_0)_{\eps>0}$ is uniformly bounded in $L^2(\R^d)$,
the families $(\psi^\eps(t))_{\eps>0}$ are uniformly bounded in $L^2(\R^d)$ for all~$t\in\R$ and we study the time evolution of their Wigner transforms defined for $(x,\xi)\in\R^{2d}$ by 
$$W^\eps(\psi^\eps(t))(x,\xi)=(2\pi)^{-d} \,\int_{\R^d} {\rm e}^{i\xi\cdot v}\psi^\eps\left(t,x-\eps \frac{v}{2}\right)\overline{\psi^\eps}\left(t,x+\eps\frac{v}{2}\right)dv, $$
and of their Wigner   measures $\mu_t$, which are the weak limits of $W^\eps(\psi^\eps(t))$ in the space of distributions (see \cite{LP}, \cite{Ge93} or the book~\cite{Zw}). 

$ $

Conical singularities naturally appear  for smooth matrix-valued potentials in the context of eigenvalue crossings. In this case, evolution laws were derived for Wigner measures in ~\cite{FG1,FL1,FG2,F4,LT,FL} and normal forms have been obtained in~\cite{C1,C2}. For single equations displaying conical singularities, transport equations were established in \cite{HLTT} in the context of acoustic waves with constant coefficients (see~\cite{Miel} for a review). Here, we are interested in a situation with variable coefficients, where the propagation phenomenon may hit the conical singularity.

\vskip 0.2cm 

\noindent Recently, Ambrosio and Figalli \cite{AF} have proposed a new approach to deal with singular potentials more general than ours. The main concept of~\cite{AF} is a regular Lagrangian flow on the space of probability measures, which allows to prove a propagation result of Wigner measures along classical trajectories {\it in average with respect to the initial data}, see~\cite{AFFG}. A related result was given
in ~\cite{FLP}, where the authors consider mixed states. On the contrary, our aim here is to consider {\it pure states}.

\vskip 0.2cm 

\noindent In this situation, we need to keep the classical point of view for singularities of the form~\aref{eq:sing}. Using the fact that classical trajectories exist and are unique we prove the propagation result for {\it every individual initial data}. In particular, we study the case of initial data with Wigner measures which concentrate on these singularities.

$ $

Wigner measures have nice geometric properties : they are measures on the cotangent space to $\R^d$, that is, on $\R^{2d}$, and they  propagate along classical trajectories as we  now recall.
Let $\mu_t$ be a Wigner measure of $(\psi^\eps(t))_{\eps>0}$ and define  
$$S=\left\{(x,\xi)\in \R^{2d},\;\;g(x)=0\right\}.$$
Since $\frac12|\xi|^2+V(x)\in{\mathcal C}^\infty(\R^{2d}\setminus S)$,
it is well-known (see \cite{GeLe93} or \cite{GMMP}), that outside~$S$ the Wigner measure satisfies the transport equation
\begin{equation}\label{eq:prop1}
\partial_t \mu_t +\nabla_x\cdot \left(\xi\, \mu_t \right)-\nabla_\xi \cdot\left(\nabla V(x) \mu_t\right)=0
\;\;{\rm in }\;\;{\mathcal D}'(S^c).
\end{equation}
The classical trajectories associated with~\aref{eq:1} are the Hamiltonian trajectories of the function~$\frac{1}{2}|\xi|^2+V(x)$, i.e. the solution curves of the ODE system
\begin{equation}\label{eq:traj}
\left\{\begin{array}{lcl}
\dot x_t(x_0,\xi_0) & = &\xi_t(x_0,\xi_0),\\
\dot \xi_t(x_0,\xi_0) & = & -\nabla V\left(x_t(x_0,\xi_0)\right),
\end{array}\right.
\end{equation}
subject to the initial conditions
$$
x_{|t=0}=x_0,\qquad \xi_{|t=0}=\xi_0.
$$ 
For $(x_0,\xi_0)\notin S$, the smoothness of $V$ near $x_0$
  implies the  existence and the uniqueness of a local solution of~\aref{eq:traj}. 
 We denote by $\Phi^t$ the flow induced by these  trajectories
 $$\Phi^t(x,\xi)=\left(x_t(x,\xi),\xi_t(x,\xi)\right),$$
 for points $(x,\xi) \notin S$ and $t$ small enough so that $\Phi_t(x,\xi)\notin S$.
  Then, equation (\ref{eq:prop1}) says, that outside~$S$ the  measure $\mu_t$  propagates along the classical trajectories as long as they do not hit~$S$. The transport equation~\aref{eq:prop1} comes from the analysis of 
 the Wigner transform  and from an Egorov type theorem   (see \cite{GMMP} or~\cite{Zw}): for $a\in{\mathcal C}_0^\infty(\R^{2d})$ and $t\in\R$ such that the support of $a\circ\Phi^{-s}$ does not intersect $S$ for all~$s\in[0,t]$, then 
\begin{equation}\label{egorof}
\langle a,W^\eps(\psi^\eps(t))\rangle\,=\,\langle a,W^\eps(\psi^\eps_0)\circ \Phi^{-t}\rangle\,+\, o(1),
\end{equation}
where the error $o(1)$ turns out to be $O(\eps ^2)$ in this context of smooth coefficients.
In this article, we study what happens when classical trajectories attain the set  
$$S^*=S\cap\{dg(x)\xi\not=0\},$$
and we extend to these points the equations~\aref{eq:prop1} and~\aref{egorof}.

$ $

We first prove that the transport equation~\aref{eq:prop1} still holds outside $S\setminus S^*$. 
 
\begin{theorem} \label{theo}
There exists a continuous map $t\mapsto \mu_t$ such that $\mu_t$ is a semi-classical measure of the family $\left(\psi^\eps(t)\right)_{\eps>0}$. Moreover, $\mu_t( S^*)=0$ for almost every $t\in\R$ and
\begin{equation}\label{eq:prop3}
\partial_t \mu_t  +\nabla_x\cdot \left(\xi \mu_t\right)-\nabla_\xi\cdot\left(\nabla V(x) \mu_t\right)=0\;\;{in}\;\;{\mathcal D'}\left(\R\times (S\setminus S^*)^c\right).
\end{equation}
\end{theorem}

Note that in view of $\mu_t( S^*)=0$ for almost every $t\in\R$,  for all $j\in\{1,\cdots,d\}$,  $\partial_{x_j} V(x)\mu_t$ is well defined as a measure on  $\R\times(S\setminus S^*)^c$.

\begin{remark} \label{rem:lambda} More precisely, we prove in Section~\ref{sec:proof} (see Remark~\ref{rem:proof13}) that
the measures $\mu_t$ satisfy
$$
\partial_t \mu_t+\nabla_x\cdot \left(\xi \mu_t\right) -\nabla_\xi\cdot\left(\nabla V(x)\mu_t\right)=\rho
$$
where $\rho$ is a distribution supported on $\R_t\times (S\setminus S^*)$ such that 
$$\exists N\in\N,\;\exists C>0,\; 
\forall a\in{\mathcal C}_0^\infty(\R^{1+2d}_{t,x,\xi}),\;\;\left|\langle a,\rho \rangle\right|\leq C\sup_{\R_t\times(S\setminus S^*)}\sup_{|\alpha|\leq N}\left| w(x)\partial_\xi^\alpha a\right|.$$
The main tool at this stage of the proof is provided by two-microlocal Wigner measures, as introduced in ~\cite{Ni}, \cite{Fe} and~\cite{F2}.  It was already the case in~\cite{HLTT} in the context of accoustic waves with constant  coefficients and in~\cite{FG1,FL1,FG2,F4,LT} when dealing with matrix-valued potentials presenting eigenvalue crossings. 
\end{remark}

The points of $S^*$ have good properties: for every point $(x_0,\xi_0)\in S^*$,  there exists a unique classical trajectory passing through it. This fact relies on the observation that if a solution of~\aref{eq:traj} satisifies  $(x_t,\xi_t)\td_t,0 (x_0,\xi_0)\in S^*$, then 
$$
\frac{g(x_t)}{|g(x_t)|} \td_t,{0^\pm} \pm \frac{dg(x_0)\xi_0}{|dg(x_0)\xi_0|}=:\omega_0.
$$
These broken trajectories are continuous but they are not ${\mathcal C}^1$. They  allow to uniquely extend the flow $\Phi^t$ to a continuous map on open sets $\Omega$ which intersects~$S$ inside~$S^*$. This generalized flow is smooth in the variable~$\xi$ as stated in the following proposition.
 
\begin{proposition}\label{prop:trajectory}
For $(x_0,\xi_0)\in S^*$ there exists $\tau_0>0$ and a unique Lipschitz continuous map $$t\mapsto \left(x_t(x_0,\xi_0),\xi_t(x_0,\xi_0)\right),\;\; t\in[-\tau_0,\tau_0]$$
satisfying~(\ref{eq:traj}) for $t\not=0$ such that $x_0(x_0,\xi_0)=x_0$, $\xi_0(x_0,\xi_0)=\xi_0$ and 
$$\displaylines{\dot x_t(x_0,\xi_0)\td_t,{0} \xi_0\;\;,\;\;\dot \xi_t(x_0,\xi_0)\td_t,{0^\pm}- \nabla V_0(x_0)\mp\,w(x_0)^tdg(x_0)\omega_0.\cr}$$
Besides, there exists a  neighborhood $\Omega$ of $(x_0,\xi_0)$ such that $\Omega\cap S\subset  S^*$ and 
 two smooth maps $[-\tau_0,\tau_0]\times (\Omega\cap S^*): t\mapsto \Phi_\pm^t(x,\xi)$ such that 
$$\forall (x,\xi)\in \Omega\cap S^*,\;\;\left(\Phi^t(x,\xi)\right)_{\pm t\in[0,\tau_0]}=\left(\Phi^t_\pm(x,\xi)\right)_{\pm t\in[0,\tau_0]}.$$
Therefore, the flow $\Phi^t$ extends to a continuous map
$$t\mapsto \Phi^t(x,\xi),\;\;t\in[-\tau_0,\tau_0],\;\;(x,\xi)\in\Omega.$$
Moreover, for $|t|<\tau_0$ and  $\alpha\in \N^d$, the maps $(x,\xi)\mapsto \partial_\xi^\alpha\Phi^t(x,\xi)$ are  continuous maps on $\Omega$ with bounded locally integrable time derivatives $\partial_t\partial_\xi^\alpha \Phi^t(x,\xi)$.
\end{proposition}

\begin{remark}
For points of $S$ which are not in $S^*$, one may lose the uniqueness of the trajectory as  the example $V(x)=-|x|$ shows: 
 the curves $x_t=\omega \frac{t^2}{2}$ and $\xi_t=\omega t$  satisfy (\ref{eq:traj}) for all $t$ and pass through $(0,0)$ at time $t=0$ independently of the choice of the vector $\omega\in {\bf S}^{d-1}$. 
 It is likely, that these non-unique trajectories induce new phenomena: in particular, the problem could become ill-posed in terms of semi-classical measures, as suggested by the example $V(x)=-|x|^{3/2}$ proposed in~\cite{AP}. 
\end{remark}

Proposition~\ref{prop:trajectory} is proved in Section~\ref{prooflemma} (note that the existence of the broken trajectories was already proved in
Proposition~1 of~\cite{FG2}). Note also that the flow~$\Phi^t$ preserves the Liouville measure close to points of $S^*$; however, besides, as a consequence of Theorem~\ref{theo} we obtain the following Theorem.  

\begin{theorem}\label{cor} 
If the initial data $(\psi^\eps_0)_{\eps>0}$ has a unique Wigner measure $\mu_0$ and if there exists $\tau_0$ such that for $t\in[0,\tau_0]$ the trajectories $\Phi^t$ issued from points of the support of $\mu_0$ do not reach $S\setminus S^*$, then $(\psi^\eps(t))_{\eps>0}$  has a unique measure $\mu_t=(\Phi^t)_*\mu_0$ for $t\in[0,\tau_0]$.
\end{theorem} 

Theorem~\ref{cor} is proved in Section~\ref{sec:cor}.  Note that the fact that $\Phi^t(x,\xi)$ is not smooth in $(x,\xi)$ close to $S^*$ makes the proof of Theorem~\ref{cor} nontrivial.  Besides, we emphasize that Theorem~\ref{cor} holds for initial data $\mu_0$ which can see $S^*$.

\vskip 0.2cm

Let us now introduce the set ${\mathcal A}$ consisting of functions $a=a(x,\xi)$  on $\R^{2d}$ such that, for every $\alpha$ with $|\alpha|\leq d+1$, the function $\partial_\xi^\alpha a$ is continuous and 
$$\left| \partial_\xi^\alpha a(x,\xi)\right|(1+|\xi|)^{d+1}\Tend{(x,\xi)}\infty 0\ ,$$
endowed with the norm
\begin{equation}\label{eq:M(a)}
M(a):={\rm max}_{|\alpha|\leq d+1}\;{\rm sup}_{(x,\xi)}\left| \partial_\xi^\alpha a(x,\xi)\right|(1+|\xi|)^{d+1} .
\end{equation}
Notice that this space is a variant of the space introduced by Lions--Paul in \cite{LP}. Then $a\circ \Phi^t\in{\mathcal A}$ for $a\in{\mathcal A}$ is compactly supported with ${\rm supp}(a)\cap S\subset S^*$, and one can consider the action of $W^\eps(\psi^\eps_0)$ on $a\circ \Phi^t$. Since the Wigner transform is convergent for the weak star topology in the dual space of ${\mathcal A}$, 
Theorem~\ref{theo} implies a weaker version of Egorov's theorem~\aref{egorof}:
For $a\in{\mathcal C}_0^\infty(\R^{2d}) $ and  $t\in\R$ such that  the support of $a\circ \Phi^{-s}$ does not intersect $S\setminus S^*$ for all $s\in[0,t]$, we have 
$$\langle a,W^\eps(\psi^\eps(t))\rangle\,=\,\langle a,W^\eps(\psi^\eps_0)\circ \Phi^{-t}\rangle\,+\,o(1).$$
However, we are not able to estimate the convergence rate in full generality.  This issue, which is interesting for numerical purpose, will be the subject of further works.

$ $

{\bf Organization of the paper}: 
The scheme of the proof of  Theorem~\ref{theo} is explained in the next Section~\ref{sec:proof}. Then, Section~\ref{sec:step1} is devoted to the analysis of the time-continuity of the measure $\mu_t$, and the transport equation is established in Section~\ref{sec:step2}; the proof of a technical lemma is the subject of Section~\ref{sec:lem}. The  analysis of the generalized flow is made in Section~\ref{sec:flow} where we prove Proposition~\ref{prop:trajectory} and the computation of the measure $\mu_t$ stated in Theorem~\ref{cor} is done in Section~\ref{sec:cor}.


\section{Scheme of the proof of  Theorem~\ref{theo}}\label{sec:proof}

 Wigner transforms are closely related to pseudodifferential operators via the formula :
$$\forall a\in{\mathcal C}_0^\infty(\R^{2d}),\;\;\forall f\in L^2(\R^d),\;\;   
\langle a,W^\eps(f )\rangle =\left(\op_\eps(a(x,\xi))f,f\right)_{L^2(\R^{d})},$$
 where the operator $\op_\eps(a)$ is the semi-classical Weyl-quantized pseudodifferential operator of symbol $a$ defined by : $\forall f\in L^2(\R^{d})$,
 \begin{equation}\label{def:pseudo}
 \op_\eps(a(x,\xi))f(x)= (2\pi)^{-d} \int_{\R^{2d} }a\left(\frac{x+x'}{2},\eps\,\xi\right) {\rm e}^{i\xi\cdot (x-x')} f(x') dx'\,d\xi,
 \end{equation}
see \cite{Zw} for example. 
Besides, by a simple adaptation of Lemma 1.1 in \cite{GeLe93} (see also Lemma~\ref{lem:estimation} below and Remark~\ref{rem:boundedness}),  one can prove that  the operator $\op_\eps(a)$ is uniformly bounded in $L^2(\R^d)$: there exists  a constant $C>0$ such that for any   $a\in L^1_{loc}(\R^{2d})$, we have
\begin{equation}\label{L2est}
 \| \op_\eps(a)\|_{{\mathcal L}(L^2(\R^d))}\leq C\,M(a),
 \end{equation}
where $M(a)$ has been defined in (\ref{eq:M(a)}). 
  The proof of the Theorem~\ref{theo} consists in three steps. 


\subsection{First step, existence of the measure}
\label{subsec:step1}
  
  Let $T>0$, we prove the existence of a sequence~$(\eps_k)$, $\eps_k\rightarrow 0$ as $k\rightarrow +\infty$, and of a continuous map $t\mapsto \mu_t$ from~$[0,T]$ into the set of positive Radon measures such that 
  for all compactly supported $a\in{\mathcal A}$  $$\forall t\in[0,T],\;\;\left(\op_{\eps_k}(a) \psi^{\eps_k}(t)\;,\;\psi^{\eps_k}(t)\right) \td_{k}, {+\infty}\int  a(x,\xi) \d\mu_t(x,\xi) .$$
  This comes from the fact (proved in Section~\ref{sec:step1}) that there exists a constant $C>0$ such that for all $a\in{\mathcal C}_0^\infty(\R^{2d})$,
  \begin{equation}\label{claim1}
  \frac{d}{dt} \left(\op_{\eps}(a) \psi^{\eps}(t)\;,\;\psi^{\eps}(t)\right)\leq C .
  \end{equation}
  Then, by diagonal extraction, considering a dense family of ${\mathcal C}_0^\infty(\R^{2d})$ and using Ascoli's Theorem, we obtain the existence of the sequence $(\eps_k)$ and of the associated family of measures $\mu_t$. Finally, we extend the convergence to compactly supported symbols $a\in{\mathcal A}$ by approaching them by $a_n\in{\mathcal C}_0^\infty(\R^{2d})$ with $M(a-a_n)\td_n,{+\infty} 0$.
  

\subsection{Second step, the transport equation}
We derive the following equation satisfied by $\mu_t$ for $t\in[0,T]$. 
  \begin{proposition}\label{prop:equation}
  There exists a distribution  $\rho$ on $[0,T]\times S\times \R^d_\xi$ such that 
  \begin{equation}\label{eq:prop4}
  \partial_t \mu_t +\nabla_x\cdot \left(\xi\mu_t \right)-\nabla_\xi\cdot \left(\nabla V(x)   {\bf 1}_{g(x)\not=0} \mu_t \right)= \rho(t,x,\xi).
  \end{equation}
  Besides, there exists $N\in\N^*$ and $C>0$ such that for all $a\in{\mathcal C}_0^\infty([0,T]\times \R^{2d})$, 
  \begin{equation}\label{*}
  \left|\langle a(t,x,\xi)\;,\;\rho\rangle \right|\, \leq \, C\, \sup_{(t,x,\xi)\in [0,T]\times S\times \R^d} \,\sup_{|\alpha|\leq N}  \left|  w(x)\partial_\xi^\alpha a (t,x,\xi) \right|\end{equation}
  where the function $w$ is defined in~\aref{eq:sing}.
  Moreover, if $\Omega$ is an open set with $\mu_t{\bf 1}_{S\cap \Omega}=0$, then for all $a$ compactly supported on $\Omega$, $\langle a,\rho\rangle=0$.
    \end{proposition}
  Proposition~\ref{prop:equation} is proved in Section~\ref{sec:step2};  the distribution $\rho$ is defined by use of two-microlocal Wigner measures  in the spirit of~\cite{Ni}, \cite{Fe} and~\cite{F2}.
 

\subsection{Third step, the measure above the singularity.}
We now prove
  \begin{equation}\label{muonS}
  \mu_t{\bf 1}_{S^*}=0
  \end{equation}
   for almost every $t\in\R$. We  
  consider the test function $a_\delta(t,x,\xi)$ depending on the small parameter $\delta\in]0,1[$, 
  $$a_\delta(t,x,\xi)= \delta \,\Phi\left(\frac{g(x)}{\delta}\right) \theta(t) b(x,\xi) $$
  where $b\in{\mathcal C}^\infty_0(\R^{2d}\setminus\{dg(x)\xi=0\})$, $b\geq 0$, $\theta\in{\mathcal C}_0^\infty([0,T])$, $\theta\geq 0$ and $\Phi\in{\mathcal C}^\infty(\R^d)$ satisfies
  $$\exists c_0>0,\;\;\forall \xi\in{\rm Supp}\, b,\;\;\forall x\in S,\;\;\nabla \Phi(0)\cdot (dg(x)\xi) >c_0.$$
 Then, in view of~\aref{*}, testing $a_\delta$ against $\rho(t,x,\xi)$ and letting $\delta$ go to $0$, we obtain 
 \begin{equation}\label{rho1}
 \langle a_\delta\;,\;\rho\rangle \td_{\delta}, 0 0.
 \end{equation}
 On the other hand, using~\aref{eq:prop4}, we obtain
 \begin{eqnarray*}
 \langle a_\delta\;,\;\rho\rangle & = & \langle a_\delta\;,\;\partial_t \mu_t +\xi\cdot \nabla_x \mu_t -\nabla V(x) \cdot \nabla_\xi \mu_t {\bf 1}_{g(x)\not=0} \mu_t  \rangle\\
 & = & \langle(dg(x)\cdot\xi)\cdot \nabla \Phi\left(\frac{g(x)}{\delta}\right) \theta(t) b(x,\xi)\;,\; \mu_t\rangle+O(\delta)
 \end{eqnarray*}
 where we have used that $\mu_t$ is a measure.  Therefore, we obtain
 \begin{equation}\label{rho2}
 \langle a_\delta\;,\;\rho\rangle \td_\delta,0  \int_{\R^{2d+1}} (dg(x)\xi)\cdot \nabla \Phi(0) \theta(t) b(x,\xi) d\mu_t(x,\xi){\bf 1}_S dt .
  \end{equation}
  In view of~\aref{rho1} and~\aref{rho2}, we have
  $$ \int_{\R^{2d+1}} (dg(x)\xi)\cdot \nabla \Phi(0) \theta(t) b(x,\xi) d\mu_t(x,\xi){\bf 1}_S dt =0.$$
   This identity implies~\aref{muonS}.
   
\subsection{Conclusion.}   
We can now conclude the proof of Theorem~\ref{theo}.  
By~\aref{muonS} and the last point of Proposition~\ref{prop:equation},  $\langle a,\rho\rangle=0$ for all $a\in{\mathcal C}_0^\infty(\R^{2d+1})$ with  ${\rm supp}(a)\cap S\subset S^*$.
Then, \aref{eq:prop4} writes \aref{eq:prop3} outside $S\setminus S^*$, which finishes the proof.

\begin{remark}\label{rem:proof13}
Note that we have proved that the distribution $\rho$ is supported above $\R\times(S\setminus S^*)$; therefore  Remark~\ref{rem:lambda} is a consequence of this observation and of  Proposition~\ref{prop:equation}.
\end{remark}


\section{Existence of the measure}\label{sec:step1}

Let us prove~\aref{claim1}. We observe that 
\begin{equation}\label{eq:mut}
\frac{d}{dt} \left(\op_\eps(a)\psi^\eps(t)\;,\;\psi^\eps(t)\right)=\frac{i}{\eps}\left(\left[ -\frac{\eps^2}{2}\Delta +V(x)\;,\; \op_\eps(a)\right]\psi^\eps(t)\;,\;\psi^\eps(t)\right).
\end{equation}
By using integration by parts, one easily obtain
\begin{equation}\label{L00}
\frac{i}{\eps}\left[- \frac{\eps^2}{2}\Delta \;,\; \op_\eps(a)\right]=\op_\eps\left(\xi\cdot\nabla_x a\right)
\end{equation}
and this family of operators is uniformly bounded.
Set
$$L_0=\frac{i}{\eps}\left[V(x)\;,\;\op_\eps(a)\right],$$
we are going to prove that this family is also uniformly bounded in $\eps$, even though~$V$ has a singularity on $S$.  
We  use the following lemma to control the norm of the considered operators. 

\begin{lemme}\label{lem:estimation}
Consider $L_\eps$ an operator of kernel $K_{\eps}(x,y)$ of the form 
\begin{equation}\label{kernel}
K_{\eps}(x,y)= \frac{1}{(2\pi\eps)^d} k\left(\frac{x+y}{2},\frac{x-y}{\eps}\right),
\end{equation}
such that the function $k$ satisfies
\begin{equation}\label{lem:cond}
N(k):=\int_{v\in\R^d} \sup_{Y\in\R^d} \left| k(Y,v)\right| dv <+\infty .
\end{equation}
Then $L_\eps$ is uniformly bounded in ${\mathcal L}(L^2(\R^d))$ and there exists $C>0$ such that $$\| L_\eps\|_{{\mathcal L}(L^2(\R^d))}\leq\, C\, N(k).$$
\end{lemme}

\begin{remark}\label{rem:boundedness}
Note that this lemma yields the uniform boundedness of the operator $\op_\eps(a)$ for $a\in{\mathcal C}_0^\infty(\R^{2d})$ and more generally for symbols $a$ compactly supported such that $\partial_\xi ^\beta a$ is bounded and locally integrable for all $|\beta|\leq d+1$.
\end{remark}

Lemma~\ref{lem:estimation}  implies 
 the boundedness of $L_0$ on $L^2(\R^d)$.
Indeed, 
 $L_0$ has a kernel $K_{\eps}(x,y)$ of the form~\aref{kernel} 
with 
$$ k_\eps(X,v)=\frac{i}{\eps}\int \left(V(X+\eps v/2)-V(X-\eps v/2)\right)  a(X,\xi){\rm e}^{i\xi\cdot v} d\xi.$$
We write
$$V(X+\eps v/2)-V(X-\eps v/2)=\eps\, G(X,\eps v)\cdot v,$$
and the boundedness of $\nabla\left(|g(x)|\right)$ on compact subsets of $\R^d$ implies the existence of  a constant $C>0$ such that
 $|G(X,\eps v)|\leq C$. Writing, thanks to an integration by parts, 
$$k_\eps(X,v)=-\int {\rm e}^{i\xi\cdot v} \nabla_\xi a(X,\xi)\cdot G(X,\eps v) d\xi$$
and using that $a$ is smooth and compactly supported in $\xi$, we obtain (again by integration by parts)
$$\forall N\in\N,\;\;\langle v \rangle ^{2N} k(X,v)= - \int G(X,\eps v) \cdot  \nabla_\xi \langle i\nabla_\xi \rangle ^{2N}a(X,\xi){\rm e}^{i\xi\cdot v} d\xi.
$$
Therefore,  we have
\begin{equation}\label{condition}
\forall N\in \N,\;\; \exists C_N>0,\;\;\sup_{Y,v} \left(\langle v \rangle ^{2N}|k_\eps (Y,v)|\right)\leq C_N 
\end{equation} 
and the  condition~\aref{lem:cond}  is satisfied. Let us now prove Lemma~\ref{lem:estimation}.

\begin{proof} 
We observe 
\begin{eqnarray*}
  \int \sup_x\left|K_{\eps}(x,y)\right| dy & = &\frac{1}{(2\pi\eps)^d} \int\sup_x \left|  k\left(\frac{x+y}{2},\frac{x-y}{\eps}\right)\right| dy\\
& = &\frac{1}{(2\pi)^d}  \int\sup_x \left| k\left(x-\eps v/2, v\right) \right|dv\\
& \leq & C  \int \sup_Y\left| k(Y,v) \right|dv.
\end{eqnarray*}
Similarly, 
\begin{eqnarray*}
\int\sup_y  \left| K_{\eps}(x,y)\right| dx & = &\frac{1}{(2\pi\eps)^d}  \int \sup_y \left|  k\left(\frac{x+y}{2},\frac{x-y}{\eps}\right)\right| dx\\
& = &\frac{1}{(2\pi)^d}  \int\sup_y  \left|k\left(y+\eps v/2, v\right)\right| dv\\
& \leq & C  \int\sup_Y\left| k(Y,v) \right|dv.
\end{eqnarray*}
Therefore, by Schur lemma, the condition \aref{condition} is enough to yield the boundedness of the operator~$L_\eps$.
\end{proof}

\section{The transport equation}\label{sec:step2}

Let us now prove Proposition~\ref{prop:equation}.
We choose $\eps=\eps_k$, the subsequence of Section~\ref{subsec:step1}.  
In view of \aref{eq:mut}, we need to pass to the limit in the term 
$L_0=\frac{i}{\eps}[V(x)\;,\;\op_\eps(a)]$.

\subsection{The smooth part.}\label{sec:smooth}
Let us consider the smooth part of the potential and set  
\begin{equation}\label{defL1}
L_1=\frac{i}{\eps}[V_0(x)\;,\;\op_\eps(a)]. 
\end{equation}
The kernel of $L_1$
is of the form~\aref{kernel} with 
\begin{eqnarray*}
k_{\eps}(X,v) &=& \frac{i}{\eps}\int \left(V_0(X+\eps v/2)-V_0(X-\eps v/2)\right)  a(X,\xi){\rm e}^{i\xi\cdot v} d\xi\\
& = &-  \int \nabla_\xi a(X,\xi) \cdot \nabla V_0(X) {\rm e}^{i\xi\cdot v} d\xi +\eps r_\eps (X,v)
\end{eqnarray*}
where $r_\eps$ satisfies
\begin{eqnarray*}
\langle v \rangle ^{2N} r_\eps(X,v) &  = &\langle v \rangle ^{2N} \frac{i}{\eps} \int a(X,\xi) {\rm e}^{i\xi\cdot v} \Bigl(V_0(X+\eps v/2)-V_0(X-\eps v/2)\\
& & \qquad \qquad \qquad \qquad 
-\eps \nabla_x V_0(X) \cdot v\Bigr) d\xi\\
& = &  \eps \int \langle i\nabla_\xi \rangle ^{2N} a(X,\xi) \Theta_\eps(X,v) d\xi
\end{eqnarray*}
with $|\Theta_\eps(X,v)| \leq C |v|^2$. Therefore, $|r_\eps(X,v)|\leq C \langle v \rangle ^{2N-2}$, and  by Lemma~\ref{lem:estimation}, the operator $\op_\eps(r_\eps)$ is uniformly bounded in $\eps$ by choosing $N$ large enough. 
As a conclusion,  we get
\begin{equation}\label{L1}
\left(L_1\psi^\eps(t)\;,\;\psi^\eps(t)\right) \td_{\eps}, 0\,- <\nabla_\xi a\cdot \nabla_xV _0\;,\; \mu_t>.
\end{equation}

\begin{remark}\label{rem:smooth1}
The previous argument only uses the following property of $V_0$:
\begin{equation}\label{key}
|V_0(X+v)-V_0(X)-\nabla_xV_0(X)v|\leq C|v|^2.
\end{equation}
 Therefore, it is enough to suppose that $V_0$ is differentiable and satisfies~\aref{key}. (The fact that  $\nabla_\xi a\cdot\nabla_x V_0\in{\mathcal A}$ is compactly supported allows us to pass to the limit in~\aref{L1}.)
\end{remark}

It remains to study the contribution of the singular part of the potential.  We first discuss the case where $g(x)=(x_1,\cdots,x_p)$, then we reduce to this situation by a local change of coordinates. The analysis  relies on a second microlocalisation on the singular set $S=\{x_1=\cdots =x_p=0\}$ in the spirit of~\cite{F2}: we explain this point in the next subsection.  

\subsection{Two microlocal Wigner measures} \label{sec:twomic}
These measures are used to describe the concentration of the family $\psi^\eps(t)$ above $S=\{x_1=\cdots =x_p=0\}$ (see~\cite{F2} for more details).
We set $$x'=(x_1,\cdots,x_p)\;\;{\rm and}\;\; x=(x',x'').$$ We consider two-microlocal test symbols $b(t,x,\xi,y)\in{\mathcal C}^\infty(\R^{2d+p+1})$ satisfying
\begin{itemize}
\item there exists a compact $K\subset \R^{2d+1}$ such that for all $y\in\R^p$, the function $(t,x,\xi)\mapsto b(t,x,\xi,y)$ is compactly supported in $K$,
\item there exists $R_0>0$ and $b_\infty(t,x,\xi,\omega)\in{\mathcal C}^\infty (\R^{2d+1}\times {\bf S}^{p-1}) $  such that for $|y|>R_0$, $b(t,x,\xi,y)=b_\infty\left(t,x,\xi,\frac{y}{|y|}\right)$,
\end{itemize}
and we analyze the action of the operator $\op_\eps(b(t,x,\xi,x'/\eps))$ as $\eps$ goes  to $0$.

\begin{proposition}\label{prop:twomic}
 There exists a positive Radon measure $\nu$ on $\R^{2d-p+1}_{t,x'',\xi} \times {\bf S}^{p-1}_\omega$ and a positive measure $M$ on $\R^{2(d-p)+1}_{t,x'',\xi''}$ valued in the set of trace-class operators on $L^2(\R^p_y)$ such that, up to a subsequence, 
$$\displaylines{
\left(\op_\eps(b(t,x,\xi,x'/\eps) \psi^\eps(t)\;,\;\psi^\eps(t)\right)  \td_\eps ,0   \int _{x'\not=0} b_\infty\left(t,x,\xi,\frac{x'}{|x'|}\right) d\mu_t(x,\xi) dt
\hfill\cr\hfill
+ \int b_\infty(t,(0,x''),\xi,\omega )d\nu(t, x'',\xi,\omega) +  {\rm tr} \int b^W(t,(0,x''),(D_y,\xi''),y) dM(t,x'',\xi'')
\cr}$$
where, for all $(x'',\xi'')\in\R^{2(d-p)}$, we denote by 
 $b^W(t,(0,x''),(D_y,\xi''),y)$ the operator obtained by the Weyl-quantization of the symbol $(y,\eta)\mapsto b(t,(0,x''),(\eta,\xi''),y)$. 
 \end{proposition}
 $ $
 
 This result (which is proved in~\cite{F2}) calls for several remarks. First, we point out that for any open set $\Omega\subset \R^{2d}$ the
 mass of the measure $\mu_t$ above $S\cap \Omega$  can be expressed in terms of the mass of $\nu$ and of the trace of $M$ according to 
\begin{eqnarray}
\label{mass}
\int _{S\cap\Omega}\mu_t(dx,d\xi) &=& \int _{\pi_{x'',\xi} (S\cap\Omega)\times{\bf S}^{p-1}} \nu(t,dx'',d\xi,d\omega)\\
&& + {\rm tr}\int_{\pi_{x'',\xi''}(S\cap\Omega)} M(t,dx'',d\xi'')\nonumber
\end{eqnarray}
where $\pi_{x'',\xi}$ and $\pi_{x'',\xi''}$ denotes the canonical projection $(x,\xi)\mapsto (x'',\xi)$ and $(x,\xi)\mapsto (x'',\xi'')$ respectively. As a consequence, $M$ and $\nu$ are measures absolutely continuous with respect to the Lebesgue measure $dt$.
Note also that for any test function $a(t,x'',\xi'')$, the operator $\langle a\;,\;M\rangle$ is a positive trace-class operator on $L^2(\R^p_y)$ so that ${\rm tr}\langle a,M\rangle \, \geq 0$; therefore, each term of the sum~\aref{mass} is positive. As a consequence, we have the following result:

\begin{remark}\label{rem:mass}  If $\mu_t{\bf 1}_{S\cap\Omega}=0$, then, by~\aref{mass} and because of the positivity of $\nu$ and~$M$, we obtain  $\nu=0$  and  $M=0$ above $\pi_{x'',\xi}(S\cap\Omega)$ and $\pi_{x'',\xi''}(S\cap \Omega)$, respectively.
\end{remark}

\noindent Moreover, we have the following characterization of the measures $\nu$ and $M$: 

\begin{remark}\label{characterization}
Let $b(t,x,\xi,y)$ be a two-microlocal test symbol  and $\chi\in{\mathcal C}_0^\infty (\R^p)$ a cut-off function such that 
\begin{equation}\label{eq:co}
\chi(y)=1\;\;\mbox{for}\;\; |y|\leq 1\;\;\mbox{and}\;\;\chi(y)=0\;\;\mbox{for}\;\;|y|\geq 2\;\;\mbox{with}\;\;0\leq \chi\leq 1. 
\end{equation}
Then, up to a subsequence in $\eps$, we have
$$\displaylines{
\limsup_{\delta},{0}\;\limsup_{R},{\infty}\; \lim_{\eps},{0} \left(\op_\eps\left(b\!\left(t,x,\xi,\frac{x'}{\eps}\right)\left(1-\chi\!\left(\frac{x'}{R\eps}\right)\right)
\chi\!\left(\frac{x'}{ \delta}\right)\right)\psi^\eps(t)\;,\;\psi^\eps(t)\right)\hfill\cr\hfill
=  \int b_\infty(t,(0,x''),\xi,\omega)d\nu(t,x,\xi,\omega), \cr
\limsup_{R},{\infty}\;\lim_{\eps},{0}\left(\op_\eps\left(b\left(t,x,\xi,\frac{x'}{\eps}\right)\chi\left(\frac{x'}{R\eps}\right)\right)\psi^\eps(t)\;,\;\psi^\eps(t)\right)\hfill\cr\hfill
=  {\rm tr} \int b^W(t,(0,x''),(D_y,\xi''),y)dM(t,x'',\xi'').\cr}$$
\end{remark}

\begin{remark}\label{rem:Mchar}
The family $(\Phi^\eps(t))_{\eps>0}$ with
\begin{equation}
\label{eq:Phi}
\Phi^\eps(t,y,x'')=\eps^{p/2} \psi^\eps(t,\eps y,x'')
\end{equation}
is uniformly bounded in $L^2(\R_{x''}^{d-p},{\mathcal H})$ where ${\mathcal H}=L^2(\R^p_y)$.
Besides, we have for $b$ compactly supported in all the variables,
$$\displaylines{\qquad \left(\op_\eps\left(b\left(t,x,\xi,\frac{x'}{\eps}\right)\right)\psi^\eps(t)\;,\;\psi^\eps(t)\right)_{L^2(\R^d)}\hfill\cr\hfill=
\left(\op_\eps\left(A_\eps(x'',\xi'')\right)\Phi^\eps(t)\;,\;\Phi^\eps(t)\right)_{L^2(\R^{d-p},{\mathcal H})}\qquad\cr}$$
where $A_\eps (x'',\xi'')=b^W\left(t,\eps y,x'',\xi'',D_y\right)$   where for $c(y,\eta)\in{\mathcal C}_0^\infty (\R^{2p})$, the operator $c^W(y,D_y)$ is the pseudodifferential operator of Weyl symbol $c(y,\eta)$  $$c^W(y,D_y)=\op_1(c(y,\eta)).$$ The operator $A_\eps(x'',\xi'')$ is a semiclassical symbol  valued in the set of  compact operators on ${\mathcal H}$, since $b(t,x',x'',\eta,\xi'',y)$ is compactly supported in $(y,\eta)$. Therefore, the measure $M$ is a semi-classical measure of the uniformly bounded family $(\Phi_\eps(t))_{\eps>0}$ of $ L^2(\R^{d-p},{\mathcal H})$.
\end{remark}
In the following subsection, we use these measures $\nu$ and $M$ to obtain a transport equation on the measure $\mu_t$.

\subsection{Concentration on a vector space} \label{sec:vectorspace}
We use the cut-off function $\chi$ of (\ref{eq:co}) and write  for $R>0$
\begin{equation}\label{def:L2L3}
\frac{i}{\eps}[|x'|w(x)\;,\;\op_\eps(a)]= L_2+L_3
\end{equation}
with 
$L_2= \frac{i}{\eps}\left[ |x'| w(x)\;,\;\op_\eps\left(a(x,\xi)\chi\left(\frac{x'}{\eps R}\right)\right)\right].$
We study separately the operators~$L_2$ and $L_3$.

$ $

\noindent {\bf Analysis of $L_2$}.  We observe
$$\left(L_2\psi^\eps(t)\;,\;\psi^\eps(t) \right)= \left(\widetilde {L_2}\Phi^\eps(t)\;,\;\Phi^\eps(t)\right)$$
with  
$$
\widetilde {L_2}=i\left[ |y| w(\eps y,x'')\;,\;\op_1\left(a(\eps y,\eps x'',\eta,\eps \xi'')\chi \left(\frac{y}{R}\right)\right)\right]
$$ and $\Phi^\eps(t)$ defined in (\ref{eq:Phi}). 
The operator $\widetilde {L_2}$ is a semiclassical operator of symbol 
$$A_\eps(x'',\xi'')= i\left[ |y| w(\eps y,x'')\;,\;a^W(\eps y,x'',D_y,\xi'')\chi \left(\frac{y}{R}\right)\right]$$
valued in the set of  operators on ${\mathcal H}$ (with the notations of Remark~\ref{rem:Mchar}).
Besides, if $\tilde \chi$ is a cut-off function such that $\tilde\chi=1$ on the support of $\chi$, we can write 
$A_\eps(x'',\xi'')=\tilde A_{\eps,R} +O(1/R)$ in operator norm
with 
$$\tilde A_{\eps,R}=i\left[ \tilde\chi \left(\frac{y}{R}\right)|y| w(\eps y,x'')\;,\;a^W(\eps y,x'',D_y,\xi'')\chi \left(\frac{y}{R}\right)\right],$$
which is a compact operator. 
By Remark~\ref{rem:Mchar}, 
for all test functions~$\theta$
\begin{eqnarray}\label{L2}
& \limsup_R,{+\infty}&  \limsup_{\eps}, 0\int \theta(t) \left(L_2\psi^\eps(t),\psi^\eps(t)\right)dt
\\
\nonumber
&= &{\rm tr}\left(\int \theta(t) i\Bigl[|y| w(0,x''), a^W(0,x'',D_y,\xi'')\Bigr] dM(t, x'',\xi'')\right).
 \end{eqnarray}

$ $

\noindent{\bf Analysis of $L_3$}.
The following lemma relates $L_3$ with the two-microlocal test symbols of subsection~\ref{sec:twomic}. 

\begin{lemme}\label{lem:L3} There exists $\eps_0>0$, $N_0\in\N$ and $C>0$ such that for all $a\in{\mathcal A}$, $\eps\in]0,\eps_0[$ and $R>1$,
$$\left\|L_3+ \op_\eps\!\left(\nabla_x(|x'|w(x))\cdot \partial_\xi 
a(x,\xi)\big(1-\chi\!\big({\textstyle\frac{x'}{R\eps}}\big)\big)\right) \right\|_{{\mathcal L}(L^2(\R^d))}\leq C M_{N_0}(a)\left(R^{-3}+\eps\right)$$
where
$$
\forall N\in\N, \;\; M_N(a)=\max_{|\alpha|\leq N}\,\sup_{(x,\xi)} \left| \partial_\xi^\alpha a(x,\xi)\right| (1+|\xi|)^{d+1}.$$
\end{lemme}

We postpone the proof to Section~\ref{sec:lem}. By Lemma~\ref{lem:L3}, we are left with the operator
$$\op_\eps\left(\nabla_x(|x'|w(x))\cdot \nabla_\xi 
a(x,\xi)\left(1-\chi\left(\frac{x'}{R\eps}\right)\right)\right).$$
Notice that the function 
$$(x,\xi)\mapsto \nabla_x(|x'|w(x))\cdot \nabla_\xi 
a(x,\xi)\left(1-\chi\left(\frac{x'}{R\eps}\right)\right)$$
is smooth.  We decompose this function in three parts :
$$\nabla_x(|x'|w(x))\cdot \nabla_\xi 
a(x,\xi)\left(1-\chi\left(\frac{x'}{R\eps}\right)\right)=b_1\left(x,\xi,\frac{x'}{\eps}\right) +b_2\left(x,\xi,\frac{x'}{\eps}\right)+c_{\eps,\delta}(x,\xi) $$
with 
\begin{eqnarray*}
b_1 (x,\xi,y)&=&w(x) \frac{y}{|y|}\cdot \nabla_{\xi'}a(x,\xi) \left(1-\chi\left(\frac{y}{R}\right)\right),\\
b_2(x,\xi,y) & = &  |x'| \nabla w(x) \cdot \nabla_\xi a(x,\xi) \left(1-\chi\left(\frac{y}{R}\right)\right)\left(1- \chi\left(\frac{x'}{\delta}\right)\right),\\
c_{\eps,\delta}(x,\xi) & =& |x'| \nabla w(x) \cdot \nabla_\xi a(x,\xi) \left(1-\chi\left(\frac{x'}{R\eps}\right)\right) \chi\left(\frac{x'}{\delta}\right).\end{eqnarray*}

$\bullet$ The symbols $b_1(x,\xi,y)$ and $b_2(x,\xi,y)$ are smooth  
 two microlocal  symbols. Therefore by Proposition~\ref{prop:twomic} (see also Remark~\ref{characterization}), we obtain that, up to a subsequence $\eps''_k$, for all test functions $\theta$ and all $j\in\{1,2\}$,
$$\displaylines{\qquad\limsup_\delta,0\; \limsup_R,{+\infty}\;  \limsup_{\eps_k''},  {0} \int \theta(t) \left( \op_\eps\left(b_j\left(x,\xi,\frac{x'}{\eps}\right)\right)\psi^\eps(t)\;,\;\psi^\eps(t)\right)dt\hfill\cr\hfill =
\, 
 \langle \theta(t) (b_j)_\infty(x,\xi,\omega) \;,\;\widetilde \rho \rangle \qquad \cr}$$
 with 
$$\displaylines{
 \widetilde \rho (t,x,\xi,\omega)  =    \mu_t (x,\xi){\bf 1}_{x'\not=0}\otimes\delta\left(\omega-\frac{x'}{|x'|}\right) + \delta(x')\otimes \nu(t,x'',\xi,\omega),\cr
 (b_1)_\infty(x,\xi,\omega)  =   w(x)\omega\cdot \nabla_{\xi'} a(x,\xi),\;\;
 (b_2)_\infty(x,\xi,\omega) =  |x'| \nabla w(x)\cdot \nabla_\xi a(x,\xi).\cr}$$

$ $

 $\bullet$ Let us now consider the symbol  $c_{\eps,\delta}$. The operator $\op_\eps(c_{\eps,\delta}(x,\xi))$ has a kernel of the form $(2\pi\eps)^{-d}k_\eps\left(\frac{x+y}{2},\frac{x-y}{\eps}\right)$ with 
 $$k_\eps(X,v)= \left(1-\chi\left(\frac{X'}{R\eps}\right)\right) \chi\left(\frac{X'}{\delta}\right)|X'|\int  \nabla w(X)\cdot \nabla_\xi a(X,\xi) {\rm e}^{i\xi\cdot v} d\xi.$$
 Therefore, using integration by parts in $\xi$, we obtain
 $$\forall N\in\N,\;\; \exists C_N>0,\;\;\langle v \rangle ^N \left| k_\eps(X,v)\right| \leq C_N \, \delta,$$
  which yields
  $$\limsup_{\delta},0\;\limsup_R,{+\infty} \;\limsup_\eps,0\;\left\|\op_\eps(c_{\eps,\delta})\right\|_{{\mathcal L}(L^2)}=0.$$
 
 $ $
 
\noindent Finally, we obtain
 \begin{eqnarray}\label{L3}
&\limsup_{R},{+\infty}&  \limsup_{\eps}, 0 \int \theta(t) \left(L_3\psi^\eps(t)\;,\;\psi^\eps(t)\right)dt    \\
\nonumber
&=&- \int \theta(t) \left(\nabla_x(|x'|w(x))\cdot \nabla_\xi a(x,\xi)\right) {\bf 1}_{x'\not=0} d\mu_t(x,\xi) dt \\
\nonumber & & -  \int \theta(t) w(x)\omega\cdot \nabla_{\xi'} a((0,x''),\xi)d\nu(t,x'',\xi,\omega) 
 \end{eqnarray}
 where we have used $(b_1)_\infty\left(x,\xi,\frac{x'}{|x'|}\right) + (b_2)_\infty\left(x,\xi,\frac{x'}{|x'|}\right)= \nabla_x(|x'|w(x))\cdot \nabla_\xi a(x,\xi).$
 
 $ $
 
\noindent As a conclusion, in view of~\aref{L2} and~\aref{L3}, we have 
\begin{eqnarray}
\label{eq:L2L3}
\lefteqn{\int \theta(t)\left((L_2+L_3)\psi^\eps(t) \;,\;\psi^\eps(t)\right)dt}\\\nonumber
&\td_\eps,0& - \int \theta(t) \left(\nabla_x(|x'|w(x))\cdot \nabla_\xi a(x,\xi)\right) {\bf 1}_{x'\not=0} d\mu_t(x,\xi) dt \\\nonumber
&& -  \int \theta(t) w(x)\omega\cdot \nabla_{\xi'} a((0,x''),\xi)d\nu(t,x'',\xi,\omega)\\\nonumber
&&+ {\rm tr}\left(\int  \theta(t) i\Bigl[ |y| w(0,x''),a^W(0,x'',D_y,\xi'')\Bigr] dM(t, x'',\xi'')\right).\nonumber
\end{eqnarray}

$ $

\noindent Let us now conclude the proof of Proposition~\ref{prop:equation}. 

\begin{proof}
We have 
\begin{eqnarray*}
\frac{d}{dt} \left(\op_\eps(a)\psi^\eps(t)\;,\;\psi^\eps(t)\right)&=&\frac{i}{\eps}\left(\left[ -\frac{\eps^2}{2}\Delta +V(x)\;,\; \op_\eps(a)\right]\psi^\eps(t)\;,\;\psi^\eps(t)\right)\\
&=&\left( \left( \op_\eps(\xi\cdot \nabla_x a) + L_1 + L_2 + L_3\right)\psi^\eps(t)\;,\;\psi^\eps(t)\right), 
\end{eqnarray*}
see~\aref{L00}, \aref{defL1} and \aref{def:L2L3}. 
By usual Weyl calculus (see for example \cite[Theorem 18.5.4.]{Ho}) the commutator resulting from $L_2$ can be written as 
$$
i\Bigl[|y| w(0,x''), a^W(0,x'',D_y,\xi'')\Bigr] = -\frac{y}{|y|}w(0,x'')(\partial_\eta a)^W(0,x'',D_y,\xi'') + r,
$$
where the symbol of $r$ depends on products of $w$ and $\eta$-derivatives of $a$. Passing to the limit $\eps\to0$, we obtain
$$
\partial_t\mu_t  =  -\xi\cdot\nabla_x\mu_t + {\bf 1}_{x'\not=0}\left(\nabla_x V_0\cdot\nabla_\xi\mu_t + \nabla_x (|x'|w(x)) \cdot\nabla_\xi \mu_t\right) + \rho 
$$ 
from \aref{L1} and \aref{eq:L2L3} for some distribution $\rho(t,x,\xi)$ satisfying the estimate \aref{*}. This is the transport equation~\aref{eq:prop4} in the case $S = \{x'=0\}$. The observation of Remark \ref{rem:mass} concludes the proof. 
\end{proof}

\subsection{More general submanifolds} 

We now suppose that $S$ is not necessarily a vector space. We work locally close to a point $x_0\in S$ in local coordinates $x=\varphi(z)$ with $z=(z',z'')\in\R^d=\R^p\times\R^{d-p}$ such that $z'=g(x)$. We consider $v^\eps= \psi^\eps \circ \varphi$. By Egorov's theorem (see~\cite[Lemma 1.10]{GeLe93} for example), the semi-classical measures $\mu$ and $\tilde \mu$ of $v^\eps$ and $\psi^\eps$, respectively, are linked by
$$\widetilde \mu(z,\zeta)= \mu(\varphi(z),\,^td\varphi(z)^{-1} \zeta).$$
Besides, $v^\eps$ solves locally, close to $x_0$,
$$i\eps\partial_t v^\eps= \op_\eps ({\textstyle\frac12}|\,^td\varphi(z)^{-1} \zeta|^2) v^\eps +( |z'| \widetilde w(z) +\widetilde V_0(z))v^\eps$$
where $\widetilde w$ and $\widetilde V_0$ are smooth. It is not difficult to check that the arguments of the preceding sections also apply  to this equation with a modified Laplacian. 
 We leave the details to the reader.


\section{Proof of Lemma~\ref{lem:L3}}\label{sec:lem}

We write $L_3=T_\eps \widetilde{L_3} T_\eps^*$ where $T_\eps$ is the scaling operator defined by
$$\forall f\in L^2(\R^d),\;\;T_\eps f(x)=\eps^{d/2} f(\eps x).$$
The we have 
$$\widetilde{L_3}= \frac{1}{i} \left[ \op_1 \left(a(\eps x,\xi)\left(1-\chi\left(\frac{x'}{R}\right)\right)\right)\;,\;|x'|w(\eps x)\right] .$$
The kernel of $\widetilde{L_3}$ is of the form $(2\pi)^{-d} K_\eps\left(\frac{x+y}{2},x-y\right)$ with 
$$\displaylines{\qquad
K_\eps(X,v)=\frac{1}{i} \int  {\rm e}^{i\xi\cdot v} a(\eps X,\xi)\left(1-\chi\left(\frac{X'}{R}\right)\right)\hfill\cr
\hfill\times\,\left[ \left|X'-\frac{v'}{2}\right|w\left(\eps X-\eps\frac{v}{2}\right) -\left| X'+\frac{v'}{2}\right|w\left(\eps X+\eps\frac{v}{2}\right)\right] d\xi.\cr}$$
We set
$$A_\eps(X,v):= \left|X'-\frac{v'}{2}\right|w\left(\eps X-\eps\frac{v}{2}\right) -\left| X'+\frac{v'}{2}\right|w\left(\eps X+\eps\frac{v}{2}\right)$$
and write 
\begin{eqnarray*}
A_\eps(X,v) &  = & \left(\left|X'-\frac{v'}{2}\right|- \left| X'+\frac{v'}{2}\right| \right) w(\eps X)\\
& &  -\frac{\eps}{2} v\cdot \nabla w(\eps X) \left(
\left|X'-\frac{v'}{2}\right|+ \left| X'+\frac{v'}{2}\right| \right) \\ & & 
+\frac{\eps ^2}{4}\left| X'-\frac{v'}{2}\right| \int_0^1 d^2 w\left(\eps X -s\eps \frac{v}{2}\right)[v,v] (1-s) ds \\
& & -\frac{\eps ^2}{4}\left| X'+\frac{v'}{2}\right| \int_0^1 d^2w\left(\eps X +s\eps \frac{v}{2}\right)[v,v] (1-s) ds \\
& = & A^{(1)}_\eps (X,v)+ A^{(2)}_\eps (X,v)+A^{(3)}_\eps (X,v)
\end{eqnarray*}
with
\begin{eqnarray}
\label {defA1}  A^{(1)}_\eps (X,v) & = &- \frac{X'}{|X'|}\cdot v' w(\eps X) -\eps v\cdot \nabla w(\eps X) |X'|\\\nonumber
&=& -\nabla(|X'|w(\eps X) )\cdot v ,\\
\label{defA2} A^{(2)}_\eps (X,v)
& =&  w(\eps X) \left(\left|X'-\frac{v'}{2}\right|- \left| X'+\frac{v'}{2}\right| +\frac{X'}{|X'|}\cdot v'\right),\\
\label{defA3} A^{(3)}_\eps (X,v)
& = & - \frac{\eps}{2} v\cdot \nabla w(\eps X) \left(\left|X'-\frac{v'}{2}\right|+ \left| X'+\frac{v'}{2}\right| -2|X'|\right)\\
\nonumber
& & 
+\frac{\eps^2}{4}\left| X'-\frac{v'}{2}\right| \int_0^1 d^2 w\left(\eps X -s\eps \frac{v}{2}\right)[v,v] (1-s) ds \\ 
\nonumber & & 
-\frac{\eps^2}{4}\left| X'+\frac{v'}{2}\right| \int_0^1 d^2 w\left(\eps X +s\eps \frac{v}{2}\right)[v,v] (1-s) ds 
\end{eqnarray}
For $j\in\{1,2,3\}$, we set
\begin{equation}\label{defKj}
K^{(j)}_\eps(X,v)  :=   \frac{1}{i} \int  {\rm e}^{i\xi\cdot v} a(\eps X,\xi)\left(1-\chi\left(\frac{X'}{R}\right)\right) A_\eps^{(j)}(X,v) d\xi.
\end{equation}
We denote by  $\widetilde{L_3^{(j)}}$ the operators  of kernel $(2\pi)^{-d} K^{(j)}_\eps \left(\frac{x+y}{2},x-y\right)$, so that we have 
\begin{equation}\label{decompL3}
\widetilde{L_3}= \widetilde{ L_3^{(1)}} + \widetilde{ L_3^{(2)}} + \widetilde{L_3^{(3)}}.
\end{equation}
We now study successively each of these operators.

\begin{remark}\label{rem:smooth2}
Here as for $V_0$ one can relax the ${\mathcal C}^2$ regularity:  assuming that $w$ is differentiable and satisfies~(\ref{key}), an argument similar to the one of the beginning of section~\ref{sec:step2} allows to perform the proof (see Remark~\ref{rem:smooth1}).
\end{remark}

$ $

$\bullet$  For $j=1$, we obtain
\begin{eqnarray*}
K^{(1)}_\eps(X,v) 
 & = & - \int  {\rm e}^{i\xi\cdot v} \nabla_\xi a(\eps X,\xi)\cdot\nabla_x \left(|X'|w(\eps X)\right)\left(1-\chi\left(\frac{X'}{R}\right)\right)  d\xi.
 \end{eqnarray*}
Therefore, the operator $\widetilde {L_3}^{(1)}$ is 
\begin{eqnarray}\nonumber
\widetilde {L_3^{(1)}} &  = &-\,  \op_1\left( \nabla_\xi a(\eps x,\xi)\cdot\nabla_x \left(|x'|w(\eps x)\right)\left(1-\chi\left(\frac{x'}{R}\right)\right) \right)\\ \label{eq:L31}
& = &-T_\eps^* \,\op_\eps \left( \nabla_\xi a(x,\xi)\cdot\nabla_x \left(|x'|w( x)\right)\left(1-\chi\left(\frac{x'}{R\eps}\right)\right) \right)T_\eps.
\end{eqnarray}

$ $

$\bullet $ For $j=2$,  we write 
\begin{eqnarray}\label{star13}\lefteqn{\left|X'-\frac{v'}{2}\right|- \left| X'+\frac{v'}{2}\right| +\frac{X'}{|X'|}\cdot v'=}\\\nonumber
&& X'\cdot v' \left( - \frac{2}{|X'+v'/2|+|X'-v'/2|}+\frac{1}{|X'|}\right).
\end{eqnarray}
Since
\begin{equation}\label{majoration}
\frac{1}{|X'+v'/2|+|X'-v'/2|}-\frac{1}{2 |X'|} =  \frac{2|X'| -|X'+v'/2| -|X'-v'/2|}{2|X'| \left( |X'+v'/2|+|X'-v'/2|\right)},
\end{equation}
we observe
\begin{eqnarray*}
\lefteqn{|X'|\left(\left|X'+\frac{v'}{2}\right|+\left|X'-\frac{v'}{2}\right| -2|X'|\right)}\\
&=& |X'|\left| X'+\frac{v'}{2}\right| +|X'|\left|X'-\frac{v'}{2}\right|-2|X'|^2\\
& \leq & \frac{1}{2}\left(\left|X'+\frac{v'}{2}\right|^2+\left|X'-\frac{v'}{2} \right|^2+2|X'|^2\right)-2|X'|^2,
\end{eqnarray*}
where we have used $ab\leq \frac{1}{2}(a^2+b^2)$. Expanding the terms $\left|X'\pm\frac{v'}{2}\right|^2$, we obtain
$$
|X'|\left(\left|X'+\frac{v'}{2}\right|+\left|X'-\frac{v'}{2}\right| -2|X'|\right)  \leq  
\frac{1}{2} \left( \frac{|v'|^2}{2}+4|X'|^2\right)-2|X'|^2= \frac{|v'|^2}{4}.$$
Plugging the latter inequality in (\ref{majoration}) and \aref{star13}, we obtain
$$\left|X'-\frac{v'}{2}\right|- \left| X'+\frac{v'}{2}\right| +\frac{X'}{|X'|}\cdot v'= \frac{X'}{|X'|}\cdot v' G(X',v')$$
with 
\begin{equation}\label{estimationG}
\left|G(X',v')\right|\leq C \, |X'|^{-3}\,|v'|^{2},
\end{equation}
where we have used that, by the triangle inequality,
$$\left| X'+\frac{v'}{2} \right|+\left|X'-\frac{v'}{2}\right|-2|X'|\geq 0$$
Therefore, by \aref{defA2}, 
$$A_\eps^{(2)}(X,v) = w(\eps X)\frac{X'}{|X'|}\cdot v' G(X',v'),$$
and integrating by parts, we have  
\begin{eqnarray*}
\lefteqn{\langle v \rangle ^{2N} K_\eps^{(2)}(X,v)=}\\
&&   \int G(X', v')\langle i\nabla_\xi \rangle ^{2N} w(\eps X) \frac{X'}{|X'|}\cdot \nabla_{\xi'} a(\eps X,\xi) \left(1-\chi\left(\frac{X'}{R}\right) \right){\rm e}^{i\xi\cdot v} d\xi .
\end{eqnarray*}
Using~\aref{estimationG} and the fact that $a$ is a smooth compactly supported function of $\xi$, we obtain that for all $N\in\N$, there exists a constant $C_N$ such that 
$$\sup_{X,v} \left| \langle v \rangle ^{2N}K_\eps^{(2)}(X,v)\right| \leq C_N \, R^{-3}.$$
We then conclude by Lemma~\ref{lem:estimation} that there exists $P\in\N$ such that 
\begin{equation}\label{eq:L32}
 \left\|\widetilde{ L_3 ^{(2)}} \right\|_{{\mathcal L}(L^2(\R^d))}\leq C M_P(a) R^{-3}.
 \end{equation}

$ $

$\bullet$ For $j=3$, we transform $A_\eps^{(3)}(X,v)$. 
We write
\begin{eqnarray*}
\lefteqn{\frac12\left(\left| X'-\frac{v'}{2}\right| +\left| X'+\frac{v'}{2}\right| -2|X'|\right)}\\
&=& \frac{v'}{4}\cdot \left( \frac{2X'+v'/2}{|X'+v'/2|+|X'|} -\frac{2X' -v'/2}{|X'-v'/2|+|X'|}\right) = v'\cdot \tilde G(X',v')
\end{eqnarray*}
with $\tilde G(X',v')$ a bounded function. Therefore,
$$A_\eps^{(3)}(X,v) = -\eps \nabla w(\eps X)\cdot v\,\tilde G(X',v')\cdot v'+\eps R_\eps(X,v)$$
with  $|R_\eps(X,v)| \leq C \langle v \rangle ^3$ for some constant $C>0$, if $\eps X$ is in the support of $a$. Finally, we obtain 
$$\displaylines{\qquad 
K_\eps^{(3)}(X,v)= - \eps \left(1-\chi\left(\frac{X'}{R}\right)\right) \int \Bigl[\nabla w(\eps X)\cdot \nabla_\xi \tilde G(X',v')\cdot \nabla_{\xi'}a(\eps X,\xi)
\hfill\cr\hfill -\eps R_\eps(X,v)a(\eps X,\xi) \Bigr]  {\rm e}^{i\xi\cdot v} d\xi,
\qquad\cr}$$
whence, by integration by parts, for all $N\in\N$
$$
\langle v \rangle ^{2N}|K_\eps^{(3)}(X,v)|  \leq  C_N\, \eps\,\sup_{k\leq 2N+4} \left| \int {\rm e}^{iv\cdot \xi } \langle i\nabla_\xi \rangle ^{k}a(\eps X,\xi ) d\xi \right| \leq C M_P(a)\eps  
$$
for some $P\in\N$. We then conclude by Lemma~\ref{lem:estimation} 
\begin{equation}\label{eq:L33}
 \|\widetilde{ L_3 ^{(3)}} \|_{{\mathcal L}(L^2(\R^d))}\leq {C} M_P(a)\,\eps.
 \end{equation}

$ $

We can now conclude the proof of Lemma~\ref{lem:L3} which comes from~\aref{decompL3} and from equations \aref{eq:L31}, \aref{eq:L32} and \aref{eq:L33}.

\section{The generalized flow}\label{sec:flow}

In the next two subsections, we  prove Proposition~\ref{prop:trajectory} in two steps: we first prove existence and uniqueness of the trajectories, then we focus on their regularity.

\subsection{Existence and uniqueness of the trajectories}\label{prooflemma}

We work with initial data $(x_0,\xi_0)\in S^*$ and  prove local existence and uniqueness of a Lipschitz map $t\mapsto (x_t,\xi_t)$ satisfying 
$\dot{x_t}=\xi_t$, $\dot{\xi_t}=-\nabla V(x_t)$ for $t\neq 0$. Then, we have
\begin{equation}
\label{eq:g} 
\frac{1}{t} g(x_t)\td_t,{0^\pm} dg(x)\xi,
\end{equation}
where $dg(x)=\left(\partial_jg_i(x)\right)_{i,j}$ is a $p\times d$ matrix and $\xi_0$ is thought as a column (a $d\times 1$ matrix); similarly, $g(x)$ is a $p\times 1$ matrix. Therefore, we have
\begin{eqnarray*}
\dot{\xi_t} &=& -\nabla V_0(x_t)-|g(x_t)|\nabla w(x_t)-w(x_t) \,^tdg(x_t)\frac{g(x_t)}{|g(x_t)|}\\
& \td_t,{0^\pm} & -\nabla V_0(x)\mp w(x) \,^tdg(x)\frac{dg(x)\xi}{|dg(x)\xi|}.
\end{eqnarray*}

\noindent For $(x,\xi)\in S^*$, we introduce the systems
$$\left\{\begin{array}{l}
\dot x^\pm_t=\xi^\pm_t,\;\;x^\pm_0=x,\\
\dot \xi^\pm_t =-\nabla V_0(x^\pm_t)\mp{\rm sgn}(t)|g(x^\pm_t)|w(x^\pm_t)\mp{\rm sgn}(t) \,^tdg(x^\pm_t)\frac{g(x^\pm_t)}{|g(x^\pm_t)|},\;\;\xi_0^\pm=\xi.
\end{array}\right.$$
We note that, under existence condition, we have 
$$\forall (x,\xi)\in\R^{2d},\;\;\Phi^t(x,\xi)=\Phi_\pm^t(x,\xi):=(x^\pm_t,\xi^\pm_t)\;\;{\rm if}\;\; \pm t>0.$$
Let us prove the existence of a solution $\Phi_+^t$, the proof for $\Phi_-^t$ is similar. We note that 
if such a map exists, then 
$t^{-1} g(x^+_t)\td_t,{0^\pm} dg(x)\xi
$
and we set
 $$y_t= \frac{1}{t} g(x^+_t)-dg(x)\xi.$$
 Since $${\rm sgn}(t) \frac{g(x^+_t)}{|g(x_t^+)|}=  \frac{t^{-1}g(x^+_t)}{|t^{-1}g(x_t^+)|}= \frac{dg(x)\xi+y_t}{|dg(x)\xi +y_t|},$$
 we are left with the system 
 $$\left\{
 \begin{array}{l}
\frac{d}{dt} (ty_t)=dg(x_t^+)\xi_t^+ - dg(x)\xi,\;\;y_0=0\\
\dot  x^+_t=\xi_t^+,\;\;x_0^+=x\\
 \dot \xi^+_t = B(t,x_t^+,y_t),\;\;\xi_0^+=\xi
 \end{array}
 \right.$$
 where
 $$B(t,X,Y):= -\nabla V_0(X)-t|dg(x)\xi+Y| -\,^t dg(X) \frac{dg(x)\xi +  Y}{|dg(x)\xi +Y|}$$
 Note that we can write
 $$ y_t=\frac{1}{t} \int_0^t \left(dg(x_s^+)\xi_s^+-dg(x)\xi\right) ds= \int_0^1  
 \left(dg(x_{t\theta}^+)\xi_{t\theta}^+-dg(x)\xi\right) d\theta.$$
Besides, since the function
$B$ is smooth near $(t,X,Y)=(0,x,0)$ for $(x,\xi)\in S^*$, we can apply a  fixed point argument to the function 
 $$\displaylines{ {\mathcal F}_{x,\xi}:\;\;(y_t,x_t^+,\xi_t^+)\hfill\cr\hfill\mapsto \left(\int_0^1  \left(dg(x_{t\theta}^+)\xi_{t\theta}^+-dg(x)\xi\right) d\theta,
x+ \int_0^t \xi_s^+ ds,\xi+ \int_0^t B(s,x_s^+,y_s)ds\right)\cr}$$
in the set ${\mathcal B}_{t_0,\delta}$ defined for $\delta,t_0>0$ small enough by
 $$ {\mathcal B}_{t_0,\delta}=\left\{ \sup_{|t|<t_0} \left(|x_t^+-x|+|\xi_t^+-\xi|+|y_t|\right) <\delta,\;y_0=0,\;x^+_0=x,\;\xi^+_0=\xi\right\}.$$
 In this manner, we construct the smooth trajectory $t\mapsto \Phi_+^t(x,\xi)$ for $(x,\xi)\in S^*$, which defines $\Phi^t(x,\xi)= \Phi_+^t(x,\xi)$ for $t>0$. The proof is similar for $t<0$ by using $\Phi_-^t$. Note also that the smoothness of $B$ with respect to $x$ and $\xi$ implies that the fixed point of ${\mathcal F}_{x,\xi}$ depends smoothly on the parameter $(x,\xi)\in S^*$.

$ $ 

\noindent These trajectories allow to uniquely extend the flow $t\mapsto\Phi^t(x,\xi)$, $t\in[-\tau_0,\tau_0]$ for $(x,\xi)\in\Omega$ with $\Omega\cap S\subset S^*$. Besides, $\partial_t \Phi^t$ is bounded and locally integrable.

\subsection{Regularity of the trajectories}

We now prove that $(x,\xi)\mapsto \partial_\xi^\alpha \Phi^t(x,\xi)$ are continuous on $\Omega$ for every multi-index~$\alpha$ by solving the system satisfied by $\partial_\xi^\alpha\Phi^t(x,\xi)$. For this, we argue by induction on~$|\alpha|$.

$ $

Let us consider $\partial_\xi^\alpha \Phi^t$ for $|\alpha|=1$ and let us first suppose $t>0$ (the argument for $t<0$ is similar).  We denote by ${\bf 1}_j$ the vector of $\N^d$ with $1$ on the $j$-th coordinate and $0$ elsewhere. We have 
$$\frac{d}{dt}  \partial_{\xi_j}\Phi^t (x,\xi)  = M(x_t) \partial_{\xi_j} \Phi^t (x,\xi) ,\;\;\partial_{\xi_j}\Phi^0(x,\xi)=(0,{\bf 1}_j)$$
 with 
 $\displaystyle{
 M(x)=
 \begin{pmatrix} 0 & {\rm Id} \\
  -B(x) & 0 \end{pmatrix}}$
and for all $\delta x\in\R^d$,
\begin{eqnarray*}
  B(x)\delta x & = & d^2V(x)\delta x\\
  & = & d^2V_0(x) \delta x + |g(x)|d^2_x w(x)\delta x+\left(\nabla w(x)\cdot \delta  x\right) \; ^tdg(x)\frac{g(x)}{|g(x)|}\\
  & & +\,w(x) d^2g(x)\left[\delta x,\frac{g(x)}{|g(x)|}\right]\\ & & 
  +\frac{w(x)}{|g(x)|}\left(
\,^tdg(x)\,dg(x)\delta x-\left(\frac{g(x)}{|g(x)|}\cdot \left(dg(x)\delta x\right)\right)\, ^tdg(x)\frac{g(x)}{|g(x)|}\right).
  \end{eqnarray*}
  Due to (\ref{eq:g}) there exists a $d\times d$ matrix $B_0$ such that 
  $$B(x_t)=\frac{B_1(x,\xi)}{t}+B_0(x,\xi)+O(t)\;\;{\rm for}\;\;t>0$$
  with
  $$B_1(x,\xi)\delta x=\frac{w(x)}{|g(x)|}\left(
\,^tdg(x)\,dg(x)\delta x-\left(\frac{dg(x)\xi}{|dg(x)\xi|}\cdot \left(dg(x)\delta x\right)\right)\, ^t dg(x)\frac{\, dg(x)\xi}{|\, dg(x)\xi|}\right).$$
  for all $\delta x\in\R^d$.
  For solving the system, we take advantage of the fact that the initial condition is such that 
  \begin{equation}\label{advantage}
 \begin{pmatrix} 0 & 0 \\
  -B_1(x,\xi) & 0 \end{pmatrix} \partial_{\xi_j}\Phi^0(x,\xi)=0.
  \end{equation}
We set 
 $$Z(t)=\frac{1}{t} \left(\partial_{\xi_j}\Phi^t-\partial_{\xi_j}\Phi^0\right)  =\frac{1}{t}\begin{pmatrix} \partial_{\xi_j}x_t\\\partial_{\xi_j}\xi_t-{\bf 1}_j\end{pmatrix}.$$
 We have 
$$\displaylines{
 \frac{d}{dt}\left( tZ(t)\right)  =  M(x_t)\left(tZ(t)+\begin{pmatrix} 0 \\ {\bf 1}_j\end{pmatrix}\right),\cr
 t\frac{d}{dt} Z(t)+Z(t)
  +  Q_0 Z(t)  =  tP(t)Z(t)+F(t)
 \;\;{\rm 
 with }\;\;
 Q_0= \begin{pmatrix} 0 & 0 \\B_1 & 0 \end{pmatrix}\cr} $$
 and 
 $t\mapsto P(t)$ and $t\mapsto F(t)$ smooth maps for $t\geq 0$. Note that there exists a unique vector $Z(0)$  such that 
 $$Z(0)+ Q_0 Z(0)=F(0).$$
We set $\tilde Z(t)=Z(t)-Z(0)$ and we have 
 $$ t\frac{d}{dt} \tilde Z(t)+\tilde Z(t)
+Q_0 \tilde Z(t)= tP(t)\tilde Z(t) +t\tilde F(t)$$
with $t\mapsto \tilde F(t)$ smooth.  We obtain 
 $$\frac{d}{dt} \left(t{\rm e} ^{Q_0{\rm ln} t}\tilde Z(t)\right)={\rm e}^{Q_0{\rm ln}t}\left(tP(t)\tilde Z(t)+t\tilde F(t)\right),$$
 where the function $t\mapsto {\rm e}^{Q_0{\rm ln} t} $ is absolutely continuous on $\R^+$,
 whence
 \begin{eqnarray*}
  t{\rm e} ^{Q_0{\rm ln} t}\tilde Z(t) & = & \int_0^t  {\rm e}^{Q_0{\rm ln}\sigma}\left(\sigma P(\sigma)\tilde Z(\sigma)+\sigma\tilde F(\sigma)\right)d\sigma\\
  \tilde Z(t) & = &t \int_0^1  {\rm e}^{Q_0{\rm ln}\theta}\left(\theta P(t\theta)\tilde Z(t\theta)+\theta\tilde F(t\theta)\right)d\theta,
  \end{eqnarray*}
  which is solved by a fixed point argument. 
 At this first step of the induction, we have obtained that the quantity
 $$\partial_{\xi_j} \Phi^t(x,\xi)=tZ(t)+\begin{pmatrix} 0\\ {\bf 1}_j\end{pmatrix}= t(Z(0)+\tilde Z(t))+ \begin{pmatrix} 0\\ {\bf 1}_j\end{pmatrix}$$
 is a continuous map on $t\geq 0$. Besides
 $$\partial_t \partial_{\xi_j} \Phi^t(x,\xi)= \tilde Z(t)+Z(0)+t\frac{d}{dt} \tilde Z(t) $$
  have finite limits when $t$ goes to $0^+$. Arguing similarly for $t\leq 0$, we build a continuous map $t\mapsto \partial_{\xi_j}\Phi^t(x,\xi)$ with a locally integrable bounded derivative $\partial_t\partial_{\xi_j}\Phi^t(x,\xi)$.

$ $

We now proceed to the last step of the induction:  we
 suppose that the functions~$t\mapsto \partial_\xi^\beta\Phi^t(x,\xi)$ 
 are well-defined for all $\beta\in\N^d$ such that  $|\beta|\leq n$ with 
  $$\displaylines{
  \partial^\beta_\xi\Phi^t(x,\xi)=O(t)\;\;{\rm if}\;\;|\beta|>1 \;\;{\rm and}\;\;\partial_{\xi_j} \Phi^t=\left(0,{\bf 1}_j\right)+O(t),\cr}  $$
 for $t$ close to $0$. Therefore, 
 $$
  \begin{pmatrix} 0 & 0 \\
  -B_1(x,\xi) & 0 \end{pmatrix}\partial_\xi^\beta \Phi^0(x,\xi) = 0.
  $$ 
 Let us consider $\partial^\alpha_\xi\Phi^t$ for $|\alpha|=n+1$. The function $\partial_\xi^\alpha\Phi^t$ satisfies  an ODE system of the form
$$
\frac{d}{dt}\partial^\alpha _\xi \Phi^t = M(x_t) \partial^\alpha_\xi\Phi^t +F(\partial^\beta_\xi \Phi^t)
$$
where the arguments of $F$ are all associated with multi-indices $\beta$  such that $|\beta|\leq n$. 
Besides, $F(\partial^\beta_\xi \Phi^t)$ is the sum of terms of the form 
$$\partial^\gamma _{x,\xi}M(x_t)\partial^{\alpha_1}_{\xi} x_t\cdots\partial^{\alpha_p}_{\xi}x_t$$ with $\gamma\in\N^{2d},\; \alpha_j\in\N^d$ and $|\alpha_1|+\cdots+|\alpha_p|=|\gamma|+1$, $|\gamma|\leq n$.   It is easy to check that 
$$\forall \alpha\in \N^{2d},\;\; \partial^\alpha_{x,\xi} B(x_t)=O(t^{-|\alpha|-1}).$$
Since $\partial^{\alpha_j}x_t=O(t)$, we obtain
$$\partial_x^\gamma B(x_t)\partial^{\alpha_1} x_t\cdots\partial^{\alpha_p}x_t=O(t^{|\alpha_1|+\cdots+|\alpha_p|-|\gamma|-1})=O(1),\;\;\gamma\in\N^d.$$
Therefore, the map $t\mapsto F(\partial^\beta_\xi \Phi_t)$ is continuous when $t$ goes to $0$. 
We are left with a system of the same form as in the first step of the induction, since the initial data have an analogous property to \aref{advantage}, and one can argue similarly, which concludes the proof of Proposition~\ref{prop:trajectory}.


\section{Propagation of the measure}\label{sec:cor}

\subsection{Preliminaries}

Before proving Theorem~\ref{cor}, we begin with a crucial remark.

\begin{remark}\label{rem:dgx}
Because of~\aref{eq:g}, the quantity $(dg(x)\xi) \cdot g(x)$ changes of sign close to $t=0$ on trajectories passing through $S^*$ at time $t=0$: $(dg(x)\xi) \cdot g(x)>0$ on the outgoing branches and $(dg(x)\xi) \cdot g(x)<0$ on the incoming ones.
\end{remark}

Let us now prove Theorem~\ref{cor}. Note that it is enough to prove the corollary under the assumption that between times $t=0$ and~$t=\tau_0$ the trajectories $\Phi^t(x,\xi)$
 issued from points of the support of $\mu_0$ cross $S^*$ at most once. We proceed in two steps: under the assumptions of  Theorem~\ref{cor}, we first  calculate $\mu_t$ near points which are 
 not in $S^*\cup \Phi^t(S^*)$, then we deal with the general case. 
  Before starting the proof, let us introduce the following notation: given a subset $A$ of $\R ^{2d+1}$ and $t\in \R $, we set  
 $$A(t):=\{ (x,\xi )\in \R ^{2d}: (t,x,\xi )\in A \} \ .$$

\subsection{The measure away from the singularity.}\label{sec:Vinout}  In this section, we prove
\begin{equation}\label{firststep}
{\bf 1}_{(S^*\cup \Phi^t(S^*))^c} \mu_t = {\bf 1}_{(S^*\cup \Phi^t(S^*))^c} (\Phi^t)_*\mu_0.
\end{equation}
We consider
$\Omega_f$ an open subset of $\R^{2d}$  such that $\Omega_f\cap S^*=\emptyset$ and a time $t_f$,  $t_f\in]0, \tau_0]$, such that there exists $t_i\in[0,t_f[$ for which the set  $\Omega_i=\Phi^{t_i-t_f}(\Omega_f)$ satisfies $\Omega_i\cap S^*=\emptyset$. It is enough to prove that $\mu_{t_f}=(\Phi^{t_f-t_i})_*\mu_{t_i}$ on $\Omega _f$. We consider the set $M$ consisting  of the points  $\left(t,\Phi^{t_f-t}(x,\xi)\right)_{t\in[t_i,t_f]}$ for all $(x,\xi)\in\Omega_f$. We have a partition of $M$, $M=V\cup V^c$, where 
$$V=\{ (t,x,\xi)\in M, \;\; \exists (s,y,\eta)\in [t_i,t_f]\times S^*,\;\; (x,\xi)=\Phi^{t-s}(y,\eta).\}$$
   The set  $V^c$ is an open subset of $\R^{2d+1}$ which is  invariant by $\Phi^t$ and we have 
 $$\forall t\in[t_i,t_f],\;\;\mu_t{\bf 1}_{V^c(t)}= (\Phi^{t-t_i})_*(\mu_{t_i}) {\bf 1}_{V^c(t)}.$$
 In particular, we have $\mu_{t_f}=(\Phi^{t_f-t_i})_*\mu_{t_i}$ in $\Omega_f\cap V^c(t_f)$. 
We will use latter that  the measure $\mu {\bf 1}_{V^c}$ also is a solution of the transport equation~\aref{eq:prop3}.  

 $ $
 
\noindent  Let us now focus on $V$.   In view of Remark~\ref{rem:dgx}, by reducing $\Omega_f$ and $t_i$ if necessary, we can assume that the quantity $(dg(x)\xi)\cdot g(x)$ vanishes in $V$ only at points of~$S^*$. Then, the set $\widetilde S:=(\R \times S)\cap M$ --- which is a subset of $ [t_i,t_f]\times S^*$ and  a submanifold of dimension $(2d-p)$ --- separates $V$ into two sets:
\begin{itemize}
\item the incoming region $V^{in}$, where  $\left(dg(x)\xi\right)\cdot g(x)<0$, which contains trajectories entering into $\widetilde S$,
\item  the outgoing region $V^{out}$, where $\left(dg(x)\xi\right)\cdot g(x) >0$, which contains trajectories which are issued from $\widetilde S$,
\end{itemize}
and we have 
$V=\widetilde S\cup V^{in} \cup V^{out}.$
Note that by the characterization through the function $(dg(x)\xi)\cdot g(x)$, the sets $V^{out}$ and $V^{in}$ have disjoint projections on $\R^{2d}$.
Because of the links between  $\Phi^t$ and $\Phi^t_\pm$,   the sets $V^{in}$ and $V^{out}$ are submanifolds of dimension $2d-p+1$
$$\displaylines{V^{in}=\{(t,x,\xi)\in M,\;\;\exists (s,y,\eta)\in [t,t_f]\times S,\;\;(x,\xi)=\Phi_-^{t-s}(y,\eta)\},\cr
V^{out}=\{(t,x,\xi)\in M,\;\;\exists (s,y,\eta)\in [t_i,t]\times S,\;\;(x,\xi)=\Phi_+^{t-s}(y,\eta)\}.\cr}$$
Finally, note that $\Omega_f\cap V(t_f)\subset V^{out}(t_f)$ and  $\Omega_i\cap V(t_i)\subset V^{in}(t_i)$ are submanifolds of dimension $2d-p$.
Note also that the vector field 
\begin{equation}\label{def:H}
H(x,\xi)=\xi\cdot \nabla_x-\nabla V(x)\cdot \nabla_\xi,
\end{equation}  is  smooth close to points $(x,\xi)$ of  $V^{in} \cup V^{out}$ and, by the definition of $V^{in} $ and $V^{out}$, it is tangent to these submanifolds. Therefore  $H$ is a vector field of $V^{in}$ and of $V^{out}$.

$ $

\includegraphics[width=11.5cm]{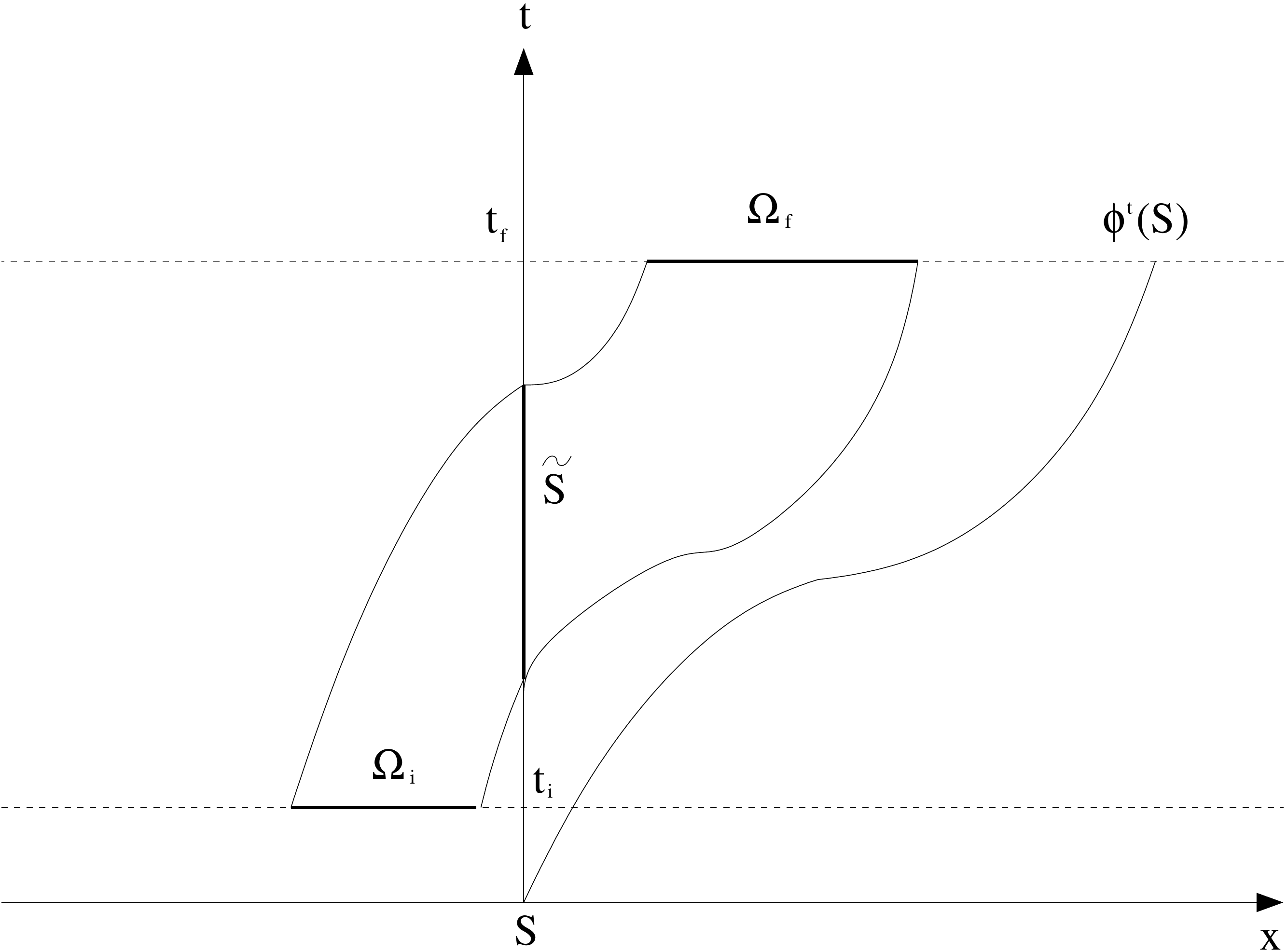}

\centerline{\bf Figure 1}

$ $

\noindent With each point $(x,\xi)$ of the projection on $T^*\R^d$ of $V$, one can associate the time $\tau(x,\xi)$ where the trajectory issued from $(x,\xi)$ passes through $S$: we have $\Phi^{\tau(x,\xi)}(x,\xi)\in S$.
If $(x,\xi)$ is  in the projection of $V^{in}$, we have $\tau(x,\xi)>0$ and  if~$(x,\xi)$ is in the projection of $V^{out}$, we have $\tau(x,\xi)<0$.

\noindent If $t_0\in [t_i,t_f]$,  we can define a map $\chi _{t_0}$ from $V(t_0)$ to $\widetilde S$ as
$$\chi_{t_0} : (x,\xi)\mapsto \left(t_0+\tau(x,\xi),\Phi^{\tau(x,\xi)}(x,\xi)\right)\in\widetilde S.$$
Note that $\chi _{t_0}$ is a homeomorphism from $V(t_0)$ onto $\widetilde S$.

$ $

\noindent Set $\mu :=\mu _t\, dt$ as a measure of $(t,x,\xi )$. We define the traces of $\mu $ on $\widetilde S$ as the measures 
$$\mu^{in}=(\chi_{t_i})_*\left(\mu_{t_i}{\bf 1}_{\Omega _i\cap V(t_i)}\right),\; \mu^{out}=(\chi_{t_f})_*\left(\mu_{t_f}{\bf 1}_{\Omega _f\cap V(t_f)}\right) .$$
Since $\mu $ satisfies the transport equation \aref{eq:prop3} on $V^{in/out}$ --- where $H$ is a smooth vector field---  and since it
does not see the set $\widetilde S$,  it is given  on $V$ by the formula
\begin{equation}\label{mu1V}
\mu {\bf 1}_V = \Phi ^{t-t_f}_*(\mu _{t_f}{\bf 1}_{t>\tau +t_f}{\bf 1}_{\Omega _f \cap V(t_f)})\, dt
+ \Phi ^{t-t_i}_*(\mu _{t_i}{\bf 1}_{t<\tau +t_i}{\bf 1}_{\Omega _i\cap V(t_i)})\, dt\;.
\end{equation}
On the other hand,   $\mu $ and $\mu {\bf 1}_{V^c}$ satisfy \aref{eq:prop3},  so  $\mu {\bf 1}_V$ does. This implies $\mu^{in}=\mu^{out}$. Indeed, the following
lemma holds.

\begin{lemme}\label{sauts}
The measure $\mu {\bf 1}_V$ satisfies the equation 
$$\partial_t (\mu {\bf 1}_V)  + \nabla_x\cdot (\xi \mu {\bf 1}_V) -\nabla_\xi \cdot (\nabla V(x) \mu  {\bf 1}_V) = {\bf 1}_{\widetilde S} (\mu^{out}-\mu^{in})\ .$$
\end{lemme}

\noindent Before proving Lemma \ref{sauts}, we observe that the relation $\mu^{out}-\mu^{in}=0$ implies
 \begin{eqnarray*}
 \mu_{t_f}{\bf 1}_{\Omega _f\cap V(t_f)}&=&(\chi_{t_f})^*\mu^{out}=(\chi_{t_f})^*\mu^{in}=(\chi_{t_f})^*(\chi_{t_{i}})_*(\mu_{t_{i}}{\bf 1}_{\Omega _i\cap V(t_i)})\\
 &=&\Phi^{t_f-t_i}_*(\mu_{t_{i}}{\bf 1}_{\Omega _i\cap V(t_i)})\ ,
 \end{eqnarray*}
 as announced.
 Let us now prove  Lemma \ref{sauts}.

\begin{proof} 
In order to compute $(\partial_t+ H) (\mu {\bf 1}_V)$, we introduce a nondecreasing function $\varphi \in C^\infty (\R )$ such that
$$\varphi (s)=0\quad {\rm for}\quad s\leq 1\ ,\ \varphi (s)=1\quad {\rm for}\quad s\ge 2\ .$$
Then
$$\displaylines{\qquad
\mu {\bf 1}_V=\lim _\delta ,{0^+}\Biggl(  \Phi ^{t-t_f}_*\left (\mu _{t_f}{\bf 1}_{\Omega _f}\varphi \left (\frac {t-t_f-\tau}{\delta} \right )\right )\, dt\hfill\cr\hfill
+\, \Phi ^{t-t_i}_*\left (\mu _{t_i}{\bf 1}_{\Omega _i}\varphi \left (\frac{t_i+\tau -t}{\delta }\right )\right )\, dt\Biggr)\ .\qquad\cr}$$
Notice that the right hand side is supported into $V^{out}\cup V^{in}$, where $H$ is smooth, so that we can make easy computations.
We obtain, in the set of distributions,
$$\displaylines{
(\partial_t+ H) (\mu {\bf 1}_V)  =\lim _\delta ,{0^+}\Biggl( \Phi ^{t-t_f}_*\left (\mu _{t_f}{\bf 1}_{\Omega _f}\frac 1\delta \varphi '\left (\frac {t-t_f-\tau}{\delta} \right )\right )\, dt\hfill\cr\hfill
-\,\Phi ^{t-t_i}_*\left (\mu _{t_i}{\bf 1}_{\Omega _i}\frac 1\delta \varphi '\left (\frac{t_i+\tau -t}{\delta }\right )\right )\, dt\Biggr)\ .\cr}$$
Therefore, given $a=a(t,x,\xi )\in C^\infty _0(M)$,
$$\displaylines{
\langle  (\partial_t+ H) (\mu {\bf 1}_V), a \rangle \hfill\cr\hfill
=\lim _\delta ,{0^+}\int _\R \int _{\Omega _f}a\left (t,\Phi ^{t-t_f}(x,\xi )\right )\frac {1}\delta \varphi '\left (\frac{t-t_f-\tau (x,\xi)}{\delta}\right )\, d\mu _{t_f}(x,\xi)\, dt\cr\hfill
-\int _\R \int _{\Omega _i}a\left (t,\Phi ^{t-t_i}(x,\xi )\right )\frac {1}\delta \varphi '\left (\frac{t_i+\tau (x,\xi )-t}{\delta}\right )\, d\mu _{t_i}(x,\xi )\, dt\ .\cr}$$
Passing to the limit in the integral as $\delta $ tends to $0^+$, we conclude
$$\langle  (\partial_t+ H) (\mu {\bf 1}_V), a \rangle= \int _{\Omega _f}a(\chi _{t_f}(x,\xi ))\, d\mu _{t_f}(x,\xi )-\int _{\Omega _i}a(\chi _{t_i}(x,\xi ))\, d\mu _{t_i}(x,\xi )$$
where we have used the definition of $\chi_t$ and the fact that $\int \varphi'(u)du=1$.
Lemma~\ref{sauts} follows by the definition of $\mu_{t_i}$ and $\mu_{t_f}$.
\end{proof}

\subsection{End of the proof of Theorem~\ref{cor}}
 We first focus on $\mu_t$ above $S^*$. Since $t\mapsto \mu_t$ and $t\mapsto \Phi^t$ are continuous, we only need to prove that for $a\in{\mathcal C}_0^\infty(\R^{2d})$ such that ${\rm supp}\left(a\right)\cap S^*=\emptyset$, we have 
\begin{equation}\label{[0,T]}
\forall T\in[0,\tau_0],\;\;
 \int_0^T\langle a\circ\Phi^{-t}, \mu_t\rangle dt = \int_0^T \langle a, \mu_0\rangle dt.
 \end{equation}
 Since $\mu ([0,T]\times S^*)=0$, we can write
$$ \int_0^T\langle a\circ\Phi^{-t}, \mu_t\rangle dt  =   \int_0^T\langle a \circ\Phi^{-t}, \mu_t {\bf 1} _{(S^*)^c}\rangle dt. $$
  Besides, since the support of $a\circ\Phi^{-t}$ does not intersect $\Phi^t(S^*)$, we have by using~\aref{firststep} 
  $$ \mu_t {\bf 1} _{(S^*)^c}=    \mu_t {\bf 1} _{(S^*)^c\cap (\Phi^t(S^*))^c}=(\Phi^t)_*\mu_0   {\bf 1} _{(S^*)^c}\;\;{\rm on}\;\;{\rm supp}\left(a\circ\Phi^{-t}\right).$$
  Therefore,  
  \begin{eqnarray*}
   \int_0^T\langle a\circ\Phi^{-t}, \mu_t\rangle dt    & = &    \int_0^T \int_{\R^{2d} } a(x,\xi)  {\bf 1} _{\Phi^{-t}((S^*)^c)}(x,\xi )\, d\mu_0(x,\xi) dt\\
  & = &\langle \int_0^T a\,  {\bf 1} _{\Phi^{-t}((S^*)^c)}  dt, \mu_0\rangle  
   \end{eqnarray*}
 where we have used the Fubini theorem.  We observe that
 $$a {\bf 1}_{\Phi^{-t}((S^*)^c)}= a - a{\bf 1} _{\Phi^{-t}(S^*)}$$
 where, for every $(x,\xi )$,   $\Phi^{t}(x,\xi)\in S^*$ for at most one value of $t$. Therefore
 $$\int_0^T a {\bf 1} _{\Phi^{-t}(S^*)} dt=0$$
 identically, and we obtain~\aref{[0,T]}. 
 
 To conclude the proof, it remains to calculate $\mu_t$ above $\Phi^t(S^*)$. For this, we work in a small neighborhood $\omega$ of a point $(x_t,\xi_t)\in\Phi^t(S^*)$. Since the flow is transverse to $S^*$, by restricting $\omega$ if necessary, we can find $\theta<0$ such that the assumptions of Theorem \ref{cor} holds on $[\theta,\tau_0]$ and such that $\Phi^{\theta-t}(\omega)\cap S^*=\emptyset$. We now argue between the times $¥\theta$ and $0$ on one hand, and between the times $\theta$ and $t$, on the other hand. The previous analysis  gives
  $$\mu_0=(\Phi^{-\theta})_*\mu_\theta\;\;{\rm on} \;\;\Phi^{-t}(\omega)\;\;{\rm
  and }\;\; \mu_t=(\Phi^{t-\theta})_*\mu_{\theta}=(\Phi^t)_*\mu_0\;\;{\rm on}\;\;\omega.$$
 This completes the proof.


\end{document}